%% file: 0_main.tex
\documentclass[12pt,a4j,twoside]{amsart}

\usepackage{docmute}
\input{setting}

\title{Convergence of energy functionals and stability of lower bounds of Ricci curvature via metric measure foliation}
\author{Daisuke Kazukawa}
\address{Mathematical Institute, Tohoku University, Sendai 980-8578, Japan} 
\email{daisuke.kazukawa.s6@dc.tohoku.ac.jp} 
\subjclass[2010]{Primary 53C23, 28A33}
\keywords{metric measure space, curvature-dimension condition, Cheeger energy functional, metric measure foliation}
\thanks{The author was supported by JSPS KAKENHI Grant Number 17J02121}
\date{\today} 

\begin{document}
\maketitle
\input{abst}
\tableofcontents
\input{sec1_intro}

\input{sec2_prelim}

\input{sec3_mmfoliation}

\input{sec4_convergence}

\input{appendix}

\input{bibliography}
\end{document}

%% file: setting.tex
\usepackage{amsmath,amssymb}
\usepackage{amsthm}
\usepackage[abbrev]{amsrefs}
\usepackage{latexsym}
\usepackage{mathrsfs}

\numberwithin{equation}{section}

\theoremstyle{definition} 
\newtheorem{thm}{Theorem}[section]
\newtheorem{prop}[thm]{Proposition}
\newtheorem{lem}[thm]{Lemma}
\newtheorem{cor}[thm]{Corollary}
\newtheorem{dfn}[thm]{Definition}
\newtheorem{ex}[thm]{Example}
\newtheorem{rem}[thm]{Remark}
\newtheorem{claim}[thm]{Claim}
\newtheorem*{ack}{Acknowledgement}

\newcommand{\CD}{\mathrm{CD}}
\newcommand{\RCD}{\mathrm{RCD}}

\newcommand{\Ent}{\mathrm{Ent}}
\newcommand{\Ch}{\mathrm{Ch}}
\newcommand{\Ric}{\mathrm{Ric}}
\newcommand{\DEnt}[1]{|D^{-} \mathrm{Ent}_{#1}|}
\newcommand{\supp}{\mathop{\mathrm{supp}}}
\newcommand{\midd}{\mathrel{} \middle| \mathrel{}}
\newcommand{\vol}{\mathrm{vol}}

\newcommand{\pr}{\mathrm{pr}}

\newcommand{\Mloc}{\mathscr{M}_{\mathrm{loc}}}
\newcommand{\Pb}{\mathscr{P}}
\newcommand{\Cbs}{C_{\mathrm{bs}}}
\newcommand{\Cb}{C_{\mathrm{b}}}
\newcommand{\resp}{\mathop{\mathrm{resp.}}}
\newcommand{\Adm}{\mathrm{Adm}}
\newcommand{\pmG}{\xrightarrow{\mathrm{pmG}}}
\newcommand{\N}{\mathbb{N}}
\newcommand{\R}{\mathbb{R}}
\newcommand{\F}{\mathcal{F}}
\newcommand{\VG}{(\mathrm{VG})}

%% file: abst.tex
\begin{abstract}
The notion of the metric measure foliation is introduced by Galaz-Garc\'ia, Kell, Mondino, and Sosa in \cite{GKMS}. They studied the relation between a metric measure space with a metric measure foliation and its quotient space. They showed that the curvature-dimension condition and the Cheeger energy functional preserve from a such space to its quotient space. Via the metric measure foliation, we investigate the convergence theory for a sequence of metric measure spaces whose dimensions are unbounded.
\end{abstract}

%% file: sec1_intro.tex
\section{Introduction}

In recent years, the geometry and analysis on metric measure spaces with Ricci curvature bounded from below are actively studied. A notion of Ricci curvature bounded from below, called the curvature-dimension condition $\CD(K, N)$, on a metric measure space has been introduced by Lott-Villani \cite{LV} and Sturm \cites{Onthe, Onthe2}. The curvature-dimension condition $\CD(K, N)$ for $K \in \mathbb{R}$ and $N \in [1, \infty]$ is defined by using the optimal transport theory and corresponds to the Ricci curvature bounded from below by $K$ and the dimension bounded from above by $N$. The class of $\CD(K, N)$ includes not only Riemannian geomtries, but also Finsler geometries. In order to isolate Riemannian from Finslerian, Ambrosio-Gigli-Savar\'e \cite{RCD} introduced the Riemannian curvature-dimension condition $\RCD(K, N)$ which is stronger than $\CD(K, N)$. 

The pmG-convergence introduced by Gigli-Mondino-Savar\'e \cite{GMS} is one of the notions of convergence of metric measure spaces. Roughly speaking, this convergence is defined by the following condition: there exists a metric space such that all metric measure spaces in a given sequence are embedded into it isometrically and the sequence of the embedded measures weakly converges. Gigli-Mondino-Savar\'e proved that the pmG-convergence is independent of the choice of embeddings and constructed the distance function metrizing the pmG-convergence on the set of all metric measure spaces. Then they proved many results for the pmG-convergence, for examples, the stability of the curvature-dimension condition, the Mosco convergences of the Cheeger energy functionals and the descending slopes of the relative entropy, the convergence of the heat flows, and the spectral convergence of the Laplacians etc.
 
The main question of our study is whether we can obtain analogous results for a sequence which does not pmG-converge. It is known that many sequences of metric measure spaces whose dimensions are unbounded do not pmG-converge. For example, the sequence of $n$-dimensional unit spheres $S^n(1)$ in $\mathbb{R}^{n+1}$, $n = 1, 2, \ldots$ with the standard Riemannian metric does not pmG-converge as $n \to \infty$ (see Corollary 5.20 and Remark 4.16 in \cite{MMG}). On the other hand, the following phenomenon occurs for sequences of $n$-dimensional spheres. For  $n$-dimensional spheres $S^n(r_n)$ of radii $r_n > 0$, we take an arbitrary point $\bar{x}_n \in S^n(r_n)$ and define a map $p_n: S^n(r_n) \to \mathbb{R}$ by
\begin{equation}\label{sphere_mmf}
p_n(x) := d_{S^n(r_n)}(x, \bar{x}_n) -\frac{\pi}{2}r_n
\end{equation}
for $x \in S^n(r_n)$, where $d_{S^n(r_n)}$ is the Riemannian distance on $S^n(r_n)$. We define a metric measure space $X_n$ for each $n$ by
\begin{equation*}
X_n := \left(\left[-\frac{\pi}{2}r_n, \frac{\pi}{2}r_n \right], |\cdot|, {p_n}_* \sigma^n \right),
\end{equation*}
where $\sigma^n$ is the normalized Riemannian volume measure on $S^n(r_n)$ and ${p_n}_* \sigma^n$ is the push-forward measure of $\sigma^n$ by $p_n$. These $X_n$ behave the following.
\begin{equation*}
\left\{ \begin{array}{ll} X_n \pmG * & \text{if } r_n/\sqrt{n} \to 0, \\
X_n \pmG (\mathbb{R}, |\cdot|, \gamma_{K^2}) & \text{if } r_n/\sqrt{n} \to K \in (0, +\infty), \\ X_n \text{ does not pmG-convergence } & \text{otherwise,} \end{array} \right.
\end{equation*}
where $*$ is a one-point metric measure space and $\gamma_{a^2}$  the 1-dimensional centered Gaussian measure on $\mathbb{R}$ with variance $a^2$. In the case that $r_n/\sqrt{n} \to K$, the Ricci curvature $\Ric_{S^n(r_n)} \equiv (n-1) / (r_n^2)$ of $S^n(r_n)$ converges to the weighted Ricci curvature $\Ric_{(\mathbb{R}, |\cdot|, \gamma_{K^2})} \equiv 1/K^2$ of the 1-dimensional Gaussian space $(\mathbb{R}, |\cdot|, \gamma_{K^2})$ of variance $K^2$  as $n \to \infty$. Moreover, for $k = 0, 1, 2, \ldots$, the $k$-th (up to multiplicity) eigenvalue $k(k + n - 1)/(r_n^2)$ of the Laplacian on $S^n(r_n)$ converges to the $k$-th eigenvalues $k/(K^2)$ of the weighted Laplacian on $(\mathbb{R}, |\cdot|, \gamma_{K^2})$  as $n \to \infty$ (see \cite{Mil}*{Subsection 2.1}). Therefore we expect our main (but still vaguely) question to be able to solve for $n$-dimensional spheres in the some sense.

Actually, the reason of the convergence of the lower bound of the Ricci curvature of these spheres has already been understood. This is that the map $p_n$ of (\ref{sphere_mmf}) induces a metric measure foliation on $S^n(r_n)$. The metric measure foliation is introduced by Galaz-Garc\'ia, Kell, Mondino, and Sosa in \cite{GKMS} and corresponds to the notion of the Riemannian submersion for metric measure spaces. The definition and other details of the metric measure foliation is described in Section 3. Galaz-Garc\'ia, Kell, Mondino, and Sosa studied the relation between a metric measure space $(X, d, m)$ with a metric measure foliation and its quotient metric measure space $(X^*, d^*, m^*)$ induced by the foliation. One of their results is that the strong curvature-dimension condition (i.e.~$\CD(K, N)$ and essentially non-branching) for $X$ implies the same condition for the quotient space $X^*$. In the observation of spheres, the space $X_n$ inherits the lower bound of the Ricci curvature from $S^n(r_n)$ and then these lower bound converges to the lower bound of the Ricci curvature of the pmG-limit space. In the smooth setting, Lott \cite{Lott} had shown that the Riemannian submersion between two weighted Riemannian manifolds preserves the lower bound of the Ricci curvature. Galaz-Garc\'ia, Kell, Mondino, and Sosa generalized this phenomenon to the framework of metric measure spaces properly. Furthermore, they showed the formula between the Cheeger energy functionals on a metric measure space $X$ with a metric measure foliation and on its quotient space $X^*$. However, their result for the Cheeger energy functional does not lead to the convergence of the eigenvalues of the Laplacian seen in the observation of spheres.

In this paper, focusing on the various convergence phenomena in the pmG-convergence, we study the metric measure foliation deeply.

The main result in this paper is the following theorem on the variational convergence of the $q$-Cheeger energy functionals $\Ch_q$. We denote by $\N$ the set of positive integers. 

\begin{thm}\label{intro_main}
Let $\{(X_n, d_n, m_n, \bar{x}_n)\}_{n \in \mathbb{N}}$ be a sequence of pointed metric measure spaces and $(Y, d, m, \bar{y})$ a pointed metric measure space and let $K \in \R$. Assume that each $X_n$ has a metric measure foliation and its quotient space $X_n^*$ satisfies the condition $\VG$ defined in Definition \ref{VG} and pmG-converges to $Y$ as $n \to \infty$. Then we have the following (1) -- (4).
\begin{enumerate}
\item If each $X_n$ satisfies $\CD(K, \infty)$, then $Y$ also satisfies $\CD(K, \infty)$.
\item Under the same assumption as in (1), $\Ch_2^{X_n}$ Mosco converges to $\Ch_2^Y$. 
\item If each $X_n$ satisfies $\RCD(K, \infty)$, then $Y$ also satisfies $\RCD(K, \infty)$.
\item Under the same assumption as in (3), $\Ch_{q_n}^{X_n}$ $\Gamma$-converges to $\Ch_q^Y$ if $\{q_n\}_{n \in \N} \subset (1, \infty)$ converges to a real number $q \in (1, \infty)$. 
\end{enumerate}
\end{thm}

\begin{rem}
\begin{enumerate}
\item The condition $\VG$ is a condition for the volume growth (see Definition \ref{VG}). This condition controls the behavior of the far measure. We also deal with a foliation with unbounded leaves not only bounded leaves.
\item In the case that $X_n$ pmG-converges to $Y$ (i.e.~each $X_n$ has the trivial foliation induced by the identity), (1) -- (3) was proved by Gigli-Mondino-Savar\'e \cite{GMS} and (4) was proved by Ambrosio-Honda \cite{AH}.
\end{enumerate}
\end{rem}

As an application of the Mosco convergence of the Cheeger energy functionals, we obtain the lower semicontinuity of the spectra of Laplacians on metric measure spaces satisfying $\RCD(K, \infty)$. The Laplacian $\Delta_X$ on a metric measure space $X$ satisfying $\RCD(K, \infty)$ is defined as the self-adjoint linear operator associated with the quadratic form $\Ch_2^X$. We denote by $\sigma(\Delta_X)$ the spectrum of $\Delta_X$.

\begin{cor}\label{intro_spec}
Under the same assumptions as in Theorem \ref{intro_main} (3), we have
\begin{equation}
\sigma(\Delta_Y) \subset \lim_{n \to \infty} \sigma(\Delta_{X_n}),
\end{equation}
that is, for any $\lambda \in \sigma(\Delta_Y)$, there exists a sequence $\lambda_n \in \sigma(\Delta_{X_n})$ convergent to $\lambda$.
\end{cor}

The following is a special case of Corollary \ref{intro_spec}. 

\begin{cor}\label{intro_spec_fc}
Let $(X, d, m)$ be a metric measure space satisfying $\RCD(K, \infty)$ for $K \in \mathbb{R}$. Assume that $X$ have a metric measure foliation and its quotient space $X^*$ satisfies $\VG$. Then we have
\begin{equation}
\sigma(\Delta_{X^*}) \subset \sigma(\Delta_X).
\end{equation}
\end{cor}

Furthermore, we obtain the $\Gamma$-convergence of the descending slopes $\DEnt{m}$ of the relative entropy $\Ent_m$. 

\begin{thm}\label{intro_main2}
Let $\{(X_n, d_n, m_n, \bar{x}_n)\}_{n \in \mathbb{N}}$ be a sequence of pointed metric measure spaces and $(Y, d, m, \bar{y})$ a pointed metric measure space and let $K \in \R$. Assume that each $X_n$ satisfies $\CD(K, \infty)$ and has a metric measure foliation and that its quotient space $X_n^*$ satisfies $\VG$ and pmG-converges to $Y$ as $n \to \infty$. Then $\DEnt{m_n}$ $\Gamma$-converges to $\DEnt{m}$.
\end{thm}

In the pmG-convergent case, Gigli-Mondino-Savar\'e proved the Mosco convergence of the descending slopes. However, since we do not know a suitable weak convergence of measures in our framework, we do not obtain the Mosco convergence. For the convergence of the heat flows, we obtain a result generalizing the result in the pmG-convergent case. The details of the descending slope of the relative entropy and the heat flow are written in the last subsection in this paper.

\begin{ack}
The author would like to thank Professor Takashi Shioya, Hiroki Nakajima, and Yuya Higashi for their comments and encouragement. He is also grateful to Professor Martin Kell for his helpful advice and for his information about Lemma \ref{Kell_lem} and its proof. He would like to thank Professor Shouhei Honda for his advice about $L^{q_n}$-convergence.
\end{ack}

%% file: sec2_prelim.tex
\section{Preliminaries}

In this section, we prepare some basic notions of the optimal transport, the Sobolev space, and the curvature-dimension condition on metric measure spaces. We use most of these notions along \cite{GMS}. As for other details, we refer to \cites{guide, V} for optimal transport, \cites{calc, G} for Sobolev space, and \cites{RCD, LV, Onthe} for curvature-dimension condition.

\subsection{Metric measure spaces and optimal transport theory}

In this paper, $(X, d)$ denotes a complete separable metric space and $m$ a locally finite Borel measure on $X$ with full support, that is, $0 < m(B_r(x)) < \infty$ for any point $x \in X$ and any real number $r > 0$. Such a triple $(X, d, m)$ is called a {\it metric measure space}, or an m.m.~space for short. In Section 4, we consider pointed metric measure spaces. We call a quadruple $(X, d, m, \bar{x})$ a {\it pointed metric measure space}, or a p.m.m.~space for short, if $(X, d, m)$ is an m.m.~space and $\bar{x} \in \supp{m}$ a base point.

We denote by $\Mloc(X)$ the set of locally finite Borel measures on $X$ and by $\Pb(X)$ the set of Borel probability measures on $X$. Further we denote by $\Cb(X)$ the set of bounded continuous functions on $X$ and by $\Cbs(X)$ the set of all functions of $\Cb(X)$ with bounded support in $X$. Then a topology of $\Mloc(X)$ is defined by the following convergence: a sequence $\{\mu_n\} \subset \Mloc(X)$ {\it converges weakly to} $\mu \in \Mloc(X)$ provided 
\begin{equation*}
\lim_{n \to \infty} \int_X \varphi(x) \, d\mu_n(x) = \int_X \varphi(x) \, d\mu(x)
\end{equation*}
for any $\varphi \in \Cbs(X)$. In the case where $\mu_n, \mu$ are finite Borel measures (for example, $\mu_n, \mu \in \Pb(X)$), the above condition is equivalent to the condition defined by $\Cb(X)$ instead of $\Cbs(X)$. 

Given two Borel probability measures $\mu_0, \mu_1 \in \Pb(X)$, we denote by $\Pi(\mu_0, \mu_1) \subset \Pb(X \times X)$ the set of transport plans between them. That means each element $\pi \in \Pi(\mu_0, \mu_1)$ satisfies ${\pr_i}_* \pi = \mu_i$ for $i = 0, 1$, where $\pr_i$ is the projection to each coordinate and ${\pr_i}_* \pi$ is the push-forward of $\pi$ by $\pr_i$.

Let $p \in [1, \infty)$ be a real number. For two probability measures $\mu, \nu \in \Pb(X)$, the {\it $L^p$-Wasserstein distance} $W_p$ between them is defined by
\begin{equation*}
W_p(\mu, \nu) := \inf_{\pi \in \Pi(\mu, \nu)} \left(\int_{X \times X} d(x, x')^p \, d\pi(x, x')\right)^{\frac{1}{p}}.
\end{equation*}
If $W_p(\mu, \nu) < + \infty$, then there exists an optimal transport plan attaining the infimum. We denote by $\Pb_p(X)$ the set of Borel probability measures on $X$ with finite $p$-th moment. Then $(\Pb_p(X), W_p)$ is a complete separable metric space and it is called the {\it $L^p$-Wasserstein space of} $X$. In some cases, we may consider the metric space $(\Pb(X), W_p)$, where the distance $W_p$ takes values in $[0, +\infty]$.

The following lemma gives a simple property of Wasserstein distance.
\begin{lem}\label{pushW_q}
Let $X, Y$ be two complete separable metric space and $p : X \to Y$ a 1-Lipschitz map and let $q \in [1, \infty)$. Then, for any $\mu_0, \mu_1 \in \Pb(X)$, we have
\begin{equation}\label{pushW_qeq}
W_q(p_\ast \mu_0, p_\ast \mu_1) \leq W_q(\mu_0, \mu_1).
\end{equation}
In other words, the map
\begin{equation}
p_\ast : \Pb(X) \ni \mu \mapsto p_\ast \mu \in \Pb(Y)
\end{equation}
is 1-Lipschitz with respect to $W_q$.
\end{lem}

\begin{proof}
We take any $\mu_0, \mu_1 \in \Pb(X)$ such that $W_q(\mu_0, \mu_1) < +\infty$. Let $\pi \in \Pb(X \times X)$ be an optimal transport plan for $W_q(\mu_0, \mu_1)$. We see that $(p \times p)_\ast \pi$ is a transport plan between $p_\ast \mu_0$ and $p_\ast \mu_1$. In fact, since $\pr_i \circ (p \times p) = p \circ \pr_i$, where $\pr_{i}$ is the projection to the $i$-th coordinate for $i = 0, 1$, we have
\begin{equation*}
{\pr_i}_\ast (p \times p)_\ast \pi = (\pr_i \circ (p \times p))_\ast \pi = (p \circ \pr_i)_\ast \pi = p_\ast {\pr_i}_\ast \pi = p_\ast \mu_i.
\end{equation*}
Therefore, 
\begin{flalign*}
W_q(p_\ast \mu_0, p_\ast \mu_1)^q \leq & \int_{Y \times Y} d_Y(y, y')^q \, d(p \times p)_\ast \pi(y, y') \\
= & \int_{X \times X} d_Y(p(x), p(x'))^q \, d\pi(x, x') \\
\leq & \int_{X \times X} d_X(x, x')^q \, d\pi(x, x') = W_q(\mu_0, \mu_1)^q.
\end{flalign*}
(\ref{pushW_qeq}) is obtained. This completes the proof.
\end{proof}

\subsection{Sobolev space on metric measure spaces}

Let $(X, d)$ be a complete separable metric space and $I \subset \mathbb{R}$  a non-trivial interval. A curve $\gamma$ on $X$ defined on $I$ means a continuous map $\gamma: I \to X$. By $C(I; X)$, we denote the space of curves on $X$ defined on $I$. We endow this space with the uniform distance and then $C(I; X)$ is a complete separable metric space.

Let $p \in [1, +\infty)$ be a real number. We define a class $AC^p(I; X)$ of curves on $X$ in the following. A curve $\gamma \in C(I; X)$ is the element of $AC^p(I; X)$ if and only if there exists $f \in L^p(I)$ satisfying
\begin{equation} \label{abs_conti}
d(\gamma(s), \gamma(t)) \leq \int_s^t f(r) \, dr
\end{equation}
for any $s, t \in I$ with $s < t$. If $p = 1$, then $\gamma$ satisfying (\ref{abs_conti}) is called an {\it absolutely continuous curve} and we write $AC(I; X)$ as $AC^1(I; X)$. For each curve $\gamma \in AC(I; X)$, it is well-known that there exists a minimal function, in the a.e.~sense, of $f$ satisfying (\ref{abs_conti}). This is called the {\it metric derivative of} $\gamma$ and is known to be provided by the following:
\begin{equation}\label{metderi}
|\dot{\gamma}|(t) := \lim_{h \to 0} \frac{d(\gamma(t+h), \gamma(t))}{|h|}
\end{equation}
for a.e.~$t \in I$ (see \cite{gradflow}*{Theorem 1.1.2}).

We define a map $\mathcal{E}_p : C(I; X) \to [0, +\infty]$ for $p > 1$ by
\begin{equation}
\mathcal{E}_p[\gamma] := \left\{ \begin{array}{ll}  \displaystyle \int_I |\dot{\gamma}|(t)^p \, dt & \text{if } \gamma \in AC^p(I; X), \\ +\infty & \text{otherwise} \end{array} \right.
\end{equation}
for any $\gamma \in C(I; X)$. The map $\mathcal{E}_p$ is lower semicontinuous and then $AC^p(I; X)$ is a Borel subset of $C(I; X)$. For $t \in I$, a continuous map $e_t : C(I; X) \to X$ is defined by $e_t(\gamma) := \gamma(t)$. 

We consider a curve $\mu : I \to \Pb(X)$ on the space $(\Pb(X), W_p)$. We often write $\mu_t$ as $\mu(t)$. Even if the distance $W_p$ takes values in $[0, +\infty]$, we can define $C(I; (\Pb(X), W_p))$ and $AC^q(I; (\Pb(X), W_p))$ for $q \in [1, +\infty)$ in the same way as above.

\begin{prop}
Let $\mu \in AC^p(I; (\Pb(X), W_p))$ and $\pi \in \Pb(C(I; X))$ satisfy ${e_t}_\ast \pi = \mu_t$ for any $t \in I$. Then, it holds that
\begin{equation}\label{lisinieq}
\int_I |\dot{\mu}|(t)^p \, dt \leq \int_{C(I; X)} \mathcal{E}_p [\gamma] \, d\pi(\gamma).
\end{equation}
\end{prop}

It is shown in \cite{Lisini} that there exists $\pi \in \Pb(C(I; X))$ satisfying equality of (\ref{lisinieq}).

\begin{prop}[\cite{Lisini}*{Corollary 1}]\label{lisinithm}
For any $\mu \in AC^p(I; (\Pb(X), W_p))$, there exists $\pi \in \Pb(C(I; X))$ such that ${e_t}_\ast \pi = \mu_t$ for any $t \in I$, and
\begin{equation}\label{lisini}
\int_I |\dot{\mu}|(t)^p \, dt = \int_{C(I; X)} \mathcal{E}_p [\gamma] \, d\pi(\gamma).
\end{equation}
\end{prop}

Let $(X, d, m)$ be an m.m.~space and let $p \in (1, \infty)$ be a real number and $q$ the conjugate exponent of $p$.

\begin{dfn}[$q$-Test plan]
We call $\pi \in \Pb(C([0,1]; X))$ a {\it $q$-test plan} provided that there exists a constant $C > 0$ such that ${e_t}_\ast \pi \leq Cm$ for any $t \in [0, 1]$, and
\begin{equation}
\int_{C([0, 1]; X)} \mathcal{E}_q[\gamma] \, d\pi(\gamma) < +\infty.
\end{equation}
\end{dfn}

\begin{dfn}[$p$-Weak upper gradient] 
Let $f : X \to \mathbb{R}$ be a Borel measurable function. We call a Borel measurable function $g : X \to [0, +\infty]$ a {\it $p$-weak upper gradient of} $f$ provided that
\begin{equation}\label{wug}
\int_{C([0,1]; X)} |f(\gamma(1)) - f(\gamma(0))| \, d\pi(\gamma) \leq \int_{C([0,1]; X)} \int_0^1 g(\gamma(t))|\dot{\gamma}|(t) \, dt d\pi(\gamma)
\end{equation}
for any $q$-test plan $\pi \in \Pb(C([0,1]; X))$. We denote by $S^p(X, d, m)$ the space of all Borel measurable functions on $X$ whose weak upper gradients belong to $L^p(X, m)$.
\end{dfn}

Given $f \in S^p(X, d, m)$, it is known that there exists a unique minimal function, in the $m$-a.e.~sense, of $p$-weak upper gradients of $f$. This is called the {\it minimal $p$-weak upper gradient of} $f$ and is denoted by $|Df|_w$, that is, for any $p$-weak upper gradient $g$, it holds that
\begin{equation}
|Df|_w(x) \leq g(x)
\end{equation}
for $m$-a.e.~$x \in X$. A more appropriate notation would be $|Df|_{w, p}$. We omit $p$ because we use only $|Df|_{w, p}$ for $f \in S^p(X, d, m)$. The details of the relation between $|Df|_{w, p}$ and $|Df|_{w, p'}$ for a function $f$ are stated in \cite{G}*{Remark 2.5}.

The {\it Sobolev space} $W^{1, p}(X, d, m)$ on an m.m.~space $(X, d, m)$ is the subspace $L^p(X, m) \cap S^p(X, d, m)$ of $L^p(X, m)$ equipped with the following norm $\| \cdot \|_{W^{1,p}}$:
\begin{equation}
\| f \|_{W^{1, p}}^p := \| f \|_{L^p}^p + \| \, |Df|_w \|_{L^p}^p.
\end{equation}

The Sobolev space $W^{1, p}(X, d, m)$ is a Banach space. However it is not a Hilbert space in general even if $p = 2$. Thus there is not always the Dirichlet form on $L^2(X, m)$ associated with the Sobolev space $W^{1,2}(X, d, m)$. Instead of the Dirichlet energy, we consider the following Cheeger energy, which is not neccesarily quadratic even if $p = 2$. We define the {\it $p$-Cheeger energy functional} $\Ch_p : L^p(X, m) \to [0, +\infty]$ by
\begin{equation}
\Ch_p(f) := \left\{ \begin{array}{ll} \displaystyle \frac{1}{p} \int_X |Df|_w (x)^p \, dm(x)  & \text{if } f \in W^{1,p}(X, d, m), \\ +\infty & \text{otherwise} \end{array} \right.
\end{equation}
for $f \in L^p(X, m)$. The functional $\Ch_p$ is lower semicontinuous and convex. 

\subsection{Curvature-dimension conditions}

Let $(X, d, m)$ be an m.m. space. The {\it relative entropy functional} $\Ent_m : \Pb(X) \to [-\infty, +\infty]$ is defined by
\begin{equation}
  \Ent_m(\mu):= \left\{ \begin{array}{ll} \displaystyle \lim_{\varepsilon \downarrow 0} \int_{\{\rho > \varepsilon \}} \rho(x) \log{\rho(x)} \, dm(x) & \text{if } \mu = \rho m, \\ +\infty & \text{otherwise} \end{array} \right.
\end{equation}
for $\mu \in \Pb(X)$. It  coincides with $\int_X \rho \log{\rho} \, dm \in [-\infty, +\infty)$ if the positive part of $\rho \log{\rho}$ is $m$-integrable, and it is equal to $+\infty$ otherwise. We denote by $D(\Ent_m)$ the set of all $\mu \in \Pb(X)$ satisfying $\Ent_m (\mu) < +\infty$. The following three properties are most important in this paper among the several basic properties of $\Ent_m$.

\begin{itemize}
\item Let $Y$ be a complete separable metric space and $p: X \to Y$ a Borel measurable map such that $p_* m \in \Mloc(Y)$. Then, for any $\mu \in \Pb(X)$, it holds that
\begin{equation}\label{pushent}
\Ent_{p_* m}(p_* \mu) \leq \Ent_{m}(\mu).
\end{equation}
\item The map $\Pb(X) \times \Pb(X) \ni (m, \mu) \mapsto \Ent_m(\mu)$ is jointly lower semicontinuous with respect to weak convergence in the two variables. 
\item Let $\mathcal{K} \subset \Mloc(X)$ and $m \in \Pb(X)$. If $\sup_{\mu \in \mathcal{K}}\Ent_m(\mu) < +\infty$, then $\mathcal{K}$ is tight.
\end{itemize}

Note that second and third properties hold only for $m \in \Pb(X)$.

The condition of Ricci curvature bounded from below on an m.m.~space $(X, d, m)$ is provided the following. 

\begin{dfn}[$\CD(K, \infty)$]\label{CD}
Let $K \in \mathbb{R}$. An m.m.~space $(X, d, m)$ satisfies the {\it curvature-dimension condition} $\CD(K, \infty)$ if for any two measures $\mu_0, \mu_1 \in \Pb_2(X) \cap D(\Ent_m)$, there exists a $W_2$-geodesic $\mu : [0,1] \ni t \mapsto \mu_t \in \Pb_2(X)$ joining $\mu_0$ and $\mu_1$ satisfying that
\begin{equation}\label{Kconvex}
\Ent_m (\mu_t) \leq (1-t) \Ent_m (\mu_0) + t \Ent_m (\mu_1) - \frac{K}{2} t(1-t) W_2(\mu_0, \mu_1)^2
\end{equation}
for any $t \in [0,1]$.
\end{dfn}

Note that a curve $\gamma : [0, 1] \to Z$ on a metric space $(Z, d)$ is called a ({\it minimal}) {\it geodesic joining} $z$ {\it and} $z'$ provided that $\gamma(0) = z$, $\gamma(1) = z'$ and
\begin{equation*}
d(\gamma(s), \gamma(t)) = |s - t| d(\gamma(0), \gamma(1))
\end{equation*}
for any $s, t \in [0, 1]$.

Let $(X, d, m)$ be an m.m.~space and assume that there exists a Lipschitz function $V : X \to [0, \infty)$  with
\begin{equation}\label{BG}
z := \int_X e^{-V(x)^2} \, dm(x) < +\infty.
\end{equation}
We define $\tilde{m} := z^{-1} e^{-V^2}m \in \Pb(X)$ and denote by $\Pb_V(X)$ the set of all $\mu \in \Pb(X)$ satisfying $\int_X V^2 \, dm < +\infty$. Since $V$ is Lipschitz, we have $\Pb_2(X) \subset \Pb_V(X)$. More generally, for any $\mu \in \Pb_V(X)$ and $\nu \in \Pb(X)$ with $W_2(\mu, \nu) < +\infty$, we have $\nu \in \Pb_V(X)$ and
\begin{equation}\label{PVW2}
\left( \int_X V^2 \, d\nu \right)^{\frac{1}{2}} \leq \mathrm{Lip}(V) W_2(\mu, \nu) + \left( \int_X V^2 \, d\mu \right)^{\frac{1}{2}},
\end{equation}
where $\mathrm{Lip}(V)$ is the Lipschitz constant of $V$.

Given $\mu \in \Pb_V(X)$, using the formula for the relative entropy
\begin{equation}\label{Entformula}
\Ent_m(\mu) = \Ent_{\tilde{m}}(\mu) - \int_X V^2 \, d\mu - \log{z},
\end{equation}
we see that $\Ent_m(\mu) > -\infty$.

\begin{dfn}[VG]\label{VG}
An m.m.~space $(X, d, m)$ satisfies the {\it volume growth condition} $\VG$ if there exist $\bar{x} \in X$ and $C > 0$ such that
\begin{equation}
\int_{X} e^{-C^2 d(x, \bar{x})^2} \, dm(x) < + \infty,
\end{equation}
that is, the Lipschitz function $V := Cd(\cdot, \bar{x})$ satisfies (\ref{BG}).
\end{dfn}

The following proposition means any $\CD(K, \infty)$ space satisfies $\VG$.

\begin{prop}[\cite{Onthe}*{Theorem 4.24}]\label{OntheVG}
Let $(X, d, m)$ be an m.m.~space satisfying $\CD(K, \infty)$ for $K \in \mathbb{R}$ and let $\bar{x} \in X$ be a fixed point. Then there exists a constant $C > 0$ such that 
\begin{equation}
m(B_r(\bar{x})) \leq C e^{(1+K_-)r^2}
\end{equation}
for every $r > 0$, where $K_- := \max{\{-K, 0 \}}$.
\end{prop}

The following lemma gives an equivalent condition for $\CD(K, \infty)$.

\begin{lem}\label{CDlem}
Let $(X, d, m)$ be an m.m.~space satisfying $\CD(K, \infty)$ for $K \in \mathbb{R}$ and let $V$ be a Lipschitz function satisfying (\ref{BG}). Then for any two measures $\mu_0, \mu_1 \in \Pb_V(X) \cap D(\Ent_m)$ with $W_2(\mu_0, \mu_1) < +\infty$, there exists a $W_2$-geodesic $\mu : [0,1] \ni t \mapsto \mu_t \in \Pb_V(X)$ joining $\mu_0$ and $\mu_1$ satisfying (\ref{Kconvex}).
\end{lem}

The proof of this lemma is described in the appendix.

\begin{dfn}[Infinitesimally Hilbertian]\label{infhilb}
An m.m.~space $(X, d, m)$ is said to be {\it infinitesimally Hilbertian} if the 2-Cheeger energy functional $\Ch_2 : L^2(X, m) \to [0, +\infty]$ is a quadratic form on $L^2(X, m)$, that is,
\begin{equation}\label{eq_infhilb}
\Ch_2(f + g) + \Ch_2(f - g) = 2 \Ch_2(f) + 2 \Ch_2(g)
\end{equation}
holds for any two functions $f, g \in L^2(X, m)$. It follows that $X$ is infinitesimally Hilbertian if and only if the Sobolev space $W^{1,2}(X, d, m)$ is a Hilbert space.
\end{dfn}

\begin{dfn}[$\RCD(K, \infty)$]\label{RCD}
Let $K \in \mathbb{R}$. An m.m.~space $(X, d, m)$ satisfies the {\it Riemannian curvature-dimension condition} $\RCD(K, \infty)$ if $X$ satisfies $\CD(K, \infty)$ and is infinitesimally Hilbertian.
\end{dfn}

%% file: sec3_mmfoliation.tex
\section{Metric measure foliation}

\subsection{Metric measure foliation}

In this subsection, we describe the metric measure foliation introduced by Galaz-Garc\'ia, Kell, Mondino and Sosa in \cite{GKMS}. We review the classical metric foliation before we explain the metric measure foliation.

\begin{dfn}[Metric foliation]
Let $(X, d)$ be a metric space and $\F$ a family of closed subsets of $X$. We call $\F$ a ({\it topological}) {\it foliation} provided that any two elements of $\F$ are disjoint to each other and $\F$ is a covering of $X$. An element $F \in \F$ is called a {\it leaf}. Furthermore, a foliation $\F$ is called a {\it metric foliation} if for any two leaves $F, F' \in \F$ and any $x \in F$,
\begin{equation}\label{mf}
d(F, F') = d(x, F').
\end{equation}
\end{dfn}

Given a metric foliation $\F$ on a metric space $(X, d)$, we consider the equivalence relation defined by $x \sim x'$ if and only if there exists $F \in \F$ such that $x, x' \in F$, and its quotient space $X^* := X/\sim$. Let $p : X \to X^*$ be the quotient map. We define a distance function $d^*$ on $X^*$ as
\begin{equation}
d^*(y, y') := d(p^{-1}(y), p^{-1}(y'))
\end{equation}
for $y, y' \in X^*$. Thanks to (\ref{mf}), the function $d^*$ becomes a distance function on $X^*$. If $(X, d)$ is complete and separable, then $(X^*, d^*)$ is also complete and separable.

We next define a submetry $f$. On a metric space $(X, d)$, we denote $B_r(x)$ the open ball centered at $x \in X$ with radius $r > 0$.

\begin{dfn}[submetry]
Let $(X, d_X), (Y, d_Y)$ be two metric spaces and $f : X \to Y$ a map between them. We call $f$ a {\it submetry} if for any $x \in X$ and $r > 0$,
\begin{equation}
f(B_r(x)) = B_r(f(x)).
\end{equation}
\end{dfn}

Note that any submetry $f$ is 1-Lipschitz and surjective.

The next lemma shows that the concepts of the submetry and the metric foliation are equivalent.

\begin{lem}[\cite{GKMS}*{Lemma 8.4}]
There is a one-to-one correspondence between metric foliations and submetries up to an isometry, that is, the following (1) and (2) hold.
\begin{enumerate}
\item Given a metric foliation $\F$ on a metric space $X$, the quotient map $p : X \to X^*$ is a submetry.
\item Given a submetry $f : X \to Y$ between two metric spaces $X$ and $Y$, the foliation $\{f^{-1}(y)\}_{y \in Y}$ is a metric foliation and there exists an isometry $i_f : Y \to X^*$ such that $i_f \circ f = p$.
\end{enumerate}
\end{lem}

In order to define the metric measure foliation, we need the disintegration obtained by the following disintegration theorem.

\begin{thm}[Disintegration theorem] \label{disint}
Let $X, Y$ be two complete separable metric spaces and $p : X \to Y$ a Borel measurable map. Then, for any Borel measure $\mu$ on $X$ satisfying $p_\ast \mu \in \Mloc(Y)$, there exists a family $\{ \mu_y \}_{y\in Y}$ of probability measures on $X$ such that
\begin{enumerate}
\item the map $Y \ni y \mapsto \mu_y(A) \in [0,1]$ is a Borel measurable function for any Borel subset $A\subset Y$,
\item $\mu_y (X \setminus p^{-1}(y)) = 0$ (i.e. $p_* \mu_y = \delta_y$) for $p_\ast \mu$-a.e. $y\in Y$,
\item for any Borel measurable function $f: X \to [-\infty, +\infty]$,
\begin{equation}\label{disinteq}
\int_X f(x) \, d \mu (x) = \int_Y \int_{p^{-1}(y)} f(x) \, d \mu_y (x) d (p_\ast \mu) (y).
\end{equation}
\end{enumerate}
Moreover, $\{ \mu_y \}_{y\in Y}$ is unique in the $p_\ast \mu$-a.e.~sense.
\end{thm}

\begin{dfn}[Disintgration] \label{disintdef}
The family $\{ \mu_y \}_{y\in Y}$ as in Theorem \ref{disint} is called the {\it disintegration} of $\mu$ for $p$.
\end{dfn}

\begin{dfn}[Metric measure foliation]\label{MMF}
Let $(X, d, m)$ be an m.m. space and $\F$ a metric foliation. We call $\F$ a {\it metric measure foliation} if $p_* m \in \Mloc(X^*)$ and there exists a Borel subset $\Omega \subset X^*$ with $p_* m(X^* \setminus \Omega) = 0$ such that  
\begin{equation}\label{mmf_eq}
W_2(\mu_y, \mu_{y'}) = d^*(y, y') = d(p^{-1}(y), p^{-1}(y'))
\end{equation}
for any $y, y' \in \Omega$, where $p : X \to X^*$ is the quotient map and $\{ \mu_y \}_{y \in Y}$ is the disintegration of $m$ for $p$.
\end{dfn}

\begin{rem}
We can weaken $p_* m \in \Mloc(X^*)$ to the condition that $p_* m$ is $\sigma$-finite because the disintegration of $m$ exists even if $p_* m$ is $\sigma$-finite. To avoid some complex situations, we deal with only locally finite measures in this papar. On the other hand, we do not assume the boundedness of the leaves which always assume in the original setting in \cite{GKMS}*{Definition 8.5}.
\end{rem}

The metric measure foliation is independent of the choice of versions of the disintegration $\{\mu_y\}_{y \in X^*}$. Moreover, since the map 
\begin{equation*}
\Omega \cap \{ y \in X^* | \, p_* \mu_y = \delta_y \} \ni y \mapsto \mu_y \in \Pb(X)
\end{equation*}
is isometric in the sense of (\ref{mmf_eq}) and $\Omega \cap \{ y \in X^* | \, p_* \mu_y = \delta_y \}$ is dense on $X^*$, this map extends to an isometric map on $X^*$. This implies that there exists a version of the disintegration $\{\mu_y\}_{y \in X^*} \subset \Pb(X)$ of $m$ such that $p_* \mu_y = \delta_y$ for all $y \in X^*$ and
\begin{equation*}
W_2(\mu_y, \mu_{y'}) = d^*(y, y')
\end{equation*}
for all two points $y, y' \in X^*$. We say that this version is {\it canonical} and often consider the canonical version of the disintegration of $m$ in the case that $\F$ is a metric measure foliation.

\begin{lem}\label{Kell_lem}
Let $(X, d_X, m_X), (Y, d_Y, m_Y)$ be two m.m.~spaces and $p : X \to Y$ a 1-Lipschitz map such that $p_* m_X = m_Y$ (so that, $p_* m_X \in \Mloc(Y)$). Assume that there exists a Borel subset $\Omega \subset Y$ with $m_Y(Y \setminus \Omega) = 0$ such that  
\begin{equation}\label{fc}
W_2(\mu_y, \mu_{y'}) = d_Y(y, y')
\end{equation}
for any $y, y' \in \Omega$. Then $p$ is a submetry. In particular, $p$ induces the metric measure foliation $\{p^{-1}(y)\}_{y \in Y}$ and $Y$ is mm-isomorphic to $X^*$. 
\end{lem}

\begin{rem}
The author studied the condition (\ref{fc}) at first. He received some important advice from Martin Kell and learned the notion of the metric measure foliation. Lemma \ref{Kell_lem} means that the condition (\ref{fc}) always induces the metric foliational structure. The proof of this lemma was given by Martin Kell.
\end{rem}

We need the following proposition for the proof of Lemma \ref{Kell_lem}.

\begin{prop}\label{optsupp}
Assume that the assumption of Lemma \ref{Kell_lem}. Let $y_0, y_1 \in Y$ be a pair of points satisfying (\ref{fc}) and $p_* \mu_{y_i} = \delta_{y_i}$ for $i=0, 1$. Then any optimal transport plan for $W_2(\mu_{y_0}, \mu_{y_1})$ is supported on
\begin{equation*}
\{(x, x') \in p^{-1}(y_0) \times p^{-1}(y_1) | \, d_Y(y_0, y_1) = d_X(x, x')\}.
\end{equation*}
\end{prop}

\begin{proof}
Let $\pi$ be an optimal transport plan for $W_2(\mu_{y_0}, \mu_{y_1})$. By the 1-Lipschitz continuity of $p$ and $\supp{\mu_{y_i}} \subset p^{-1}(y_i)$ for $i=0, 1$, we see that $\pi$ is supported on 
\begin{equation*}
\{(x, x') \in p^{-1}(y_0) \times p^{-1}(y_1) | \, d_Y(y_0, y_1) \leq d_X(x, x')\}.
\end{equation*}
If $\Gamma := \{(x, x') \in p^{-1}(y_0) \times p^{-1}(y_1) | \, d_Y(y_0, y_1) < d_X(x, x')\}$ is not $\pi$-negligible, then
\begin{equation*}
\pi(\Gamma) \cdot d_Y(y_0,y_1)^2 < \int_{\Gamma} d_X(x, x')^2 \, d\pi(x, x').
\end{equation*}
Therefore we have 
\begin{equation*}
d_Y(y_0, y_1)^2 < \int_{X \times X} d_X(x, x')^2 \, d\pi(x, x') = W_2(\mu_{y_0}, \mu_{y_1})^2,
\end{equation*}
which contradicts (\ref{fc}). This completes the proof.
\end{proof}

\begin{proof}[{\bf Proof of Lemma \ref{Kell_lem}}]
Let $\{ \mu_y \}_{y \in Y}$ be the canonical disintegration of $m_X$. Note that we obtain the canonical one using only (\ref{fc}). In particular, the property $p_* \mu_y = \delta_y$ for all $y \in Y$ implies that $p$ is surjective.

Suppose that there are $x \in X$ and $r > 0$ such that $B_r(y) \setminus p(B_r(x)) \neq \emptyset$, where $y = p(x)$. By this assmption, there exist a neighborhood $U$ of $x$ and a real number $r' < r$ such that $B_{r'}(y) \setminus p(V_{r'}) \neq \emptyset$, where
\begin{equation*}
V_{r'} := \bigcup_{\tilde{x} \in U \cap p^{-1}(y)} B_{r'}(\tilde{x}).
\end{equation*}
In fact, by $B_r(y) \setminus p(B_r(x)) \neq \emptyset$, there exist $y' \in Y$ and $r' < r$ such that $y' \in B_{r'}(y) \setminus p(B_{r} (x))$. Setting $U := B_{r-r'}(x)$ and $V_{r'}$ as above, we have $V_{r'} \cap p^{-1}(y') = \emptyset$. If there is $x' \in V_{r'} \cap p^{-1}(y')$, then there exists $\tilde{x} \in U \cap p^{-1}(y)$ such that $x' \in B_{r'}(\tilde{x})$ and
\begin{equation*}
d_X(x, x') \leq d_X(x, \tilde{x}) + d_X(\tilde{x}, x') < (r - r') + r' = r,
\end{equation*}
which implies $y' \in p(B_r(x))$. This contradicts the choice of $y'$. Thus $V_{r'} \cap p^{-1}(y') = \emptyset$, that is, $y' \in B_{r'}(y) \setminus p(V_{r'})$. 

Let $\pi_{yy'}$ be an optimal transport plan for $W_2(\mu_y, \mu_{y'})$ for some $y' \in B_{r'}(y) \setminus p(V_{r'})$. By the choice of $y'$, we see that
\begin{equation*}
d_Y(y, y') < r' \leq d_X(\tilde{x}, x')
\end{equation*}
for all $(\tilde{x}, x') \in (U \cap p^{-1}(y)) \times p^{-1}(y')$. Since $m_X$ has full support, we may assume that $\mu_y(U \cap p^{-1}(y)) > 0$. Thus we have 
\begin{align*}
&\pi_{yy'}(\{(\tilde{x}, x') \in p^{-1}(y) \times p^{-1}(y') | \, d_Y(y, y') < d_X(\tilde{x}, x')\}) \\
&\geq \pi_{yy'}((U \cap p^{-1}(y)) \times p^{-1}(y')) = \mu_y(U \cap p^{-1}(y)) > 0.
\end{align*}
This contradicts Proposition \ref{optsupp}. The proof is completed.
\end{proof}

\begin{prop}
Let $(X, d, m)$ be an m.m.~space and $\F$ a metric measure foliation. Then we have
\begin{equation}
W_q(\mu_y, \mu_{y'}) = d^*(y, y') = d(p^{-1}(y), p^{-1}(y'))
\end{equation}
for any $q \in [1, \infty)$ and any $y, y' \in X^*$ where $\{\mu_y\}_{y \in X^*}$ is the canonical disintegration.
\end{prop}

\begin{proof}
We take any $q \in [1, \infty)$ and any $y, y' \in X^*$. Let $\pi_{yy'}$ be an optimal transport plan for $W_2(\mu_y, \mu_{y'})$. By Proposition \ref{optsupp}, we see that $\pi_{yy'}$ is supported on
\begin{equation*}
\{(x, x') \in p^{-1}(y) \times p^{-1}(y') | \, d^*(y, y') = d(x, x')\},
\end{equation*}
which means that
\begin{equation*}
W_q(\mu_y, \mu_{y'})^q \leq \int_{X \times X} d(x, x')^q \, d\pi_{yy'}(x, x') = d^*(y, y')^q.
\end{equation*}
On the other hand, by Lemma \ref{pushW_q}, we see that $d^*(y, y') \leq W_q(\mu_y, \mu_{y'})$. The proof is completed.
\end{proof}

We give some examples of metric measure foliation.

\begin{ex}[Riemannian submersion]
The Riemannian submersion between two weighted Riemannian manifolds induces a metric measure foliation. The notion of the metric measure foliation is motivated from Lott's article \cite{Lott} about a relation between the weighted Ricci curvature and the Riemannian submersion. This detail is described in \cite{GKMS}.

Let $(M, g, \varphi \vol_g), (N, h, \psi \vol_h)$ be two weighted Riemannian manifolds and $\pi : M \to N$ a Riemannian submersion such that $\pi_* (\varphi \vol_g) = \psi \vol_h$. For any smooth curve $\gamma : [0,1] \to N$, we define a diffeomorphism $\rho : \pi^{-1}(\gamma(0)) \to \pi^{-1}(\gamma(1))$ between the fibers of the extremal points of $\gamma$ as the correspondence of the two extremal points of each horizontal lift of $\gamma$, that is, $\rho_\gamma(x) := \bar{\gamma}_x(1)$, where $\bar{\gamma}_x$ is the horizontal lift of $\gamma$ with $\bar{\gamma}_x(0) = x$.

Assume that $\rho_* \mu_{\gamma(0)} = \mu_{\gamma(1)}$ for any smooth curve $\gamma$ on $N$. Then the family $\{\pi^{-1}(y)\}_{y \in N}$ is a metric measure foliation on $M$.
\end{ex}

\begin{ex}[$l_q$-Product space]
The product space of m.m.~spaces is a typical example of a metric measure foliation. Let $(Y, d_Y, m_Y)$, $(Z, d_Z, m_Z)$ be two m.m.~spaces and $q \in [1, +\infty]$ an extended real number. We define the {\it $l_q$-product} $Y \times_{l_q} Z$ of $Y$ and $Z$ as the product space $Y \times Z$ equipped with the distance $d_{l_q}$ and the measure $m_Y \otimes m_Z$, where $d_{l_q}$ is defined by
\begin{equation}
d_{l_q}((y, z), (y', z')) := \left\{ \begin{array}{ll} \displaystyle (d_Y(y, y')^q + d_Z(z, z')^q)^\frac{1}{q} &  \text{if } 1 \leq q < +\infty \\ \max \{d_Y(y, y'), d_Z(z, z')\} & \text{if } q = +\infty \end{array} \right.
\end{equation}
for any two points $(y, z), (y', z') \in Y \times Z$, and $m_Y \otimes m_Z$ means the product measure of $m_Y$ and $m_Z$. 

The $l_p$-product space $Y \times_{l_q} Z$ has a metric measure foliation induced by the projection $p : Y \times_{l_q} Z \to Y$ if $m_Z$ has  the finite mass.
\end{ex}

\begin{ex}[Action of isometry group]
An m.m.~space with an isometric action by a compact group is an important example of the metric measure foliation. This is studied in \cite{GKMS} deeply.

Let $(X, d, m)$ be an m.m.~space and $G$ a compact (topological) group. Let $G \times X \ni (g, x) \mapsto gx \in X$ be an isometric action of $G$ on $X$. Then, the distance function $d_{X/G}$ on the quotient space $X/G$ is defined by
\begin{equation}
d_{X/G}([x], [x']) = \inf_{g,g' \in G} d(gx, g'x')
\end{equation}
for $[x], [x'] \in X/G$, where $[x]$ is the $G$-orbit of a point $x \in X$. $(X/G, d_{X/G})$ is a complete separable metric space. Let $p : X \ni x \mapsto [x] \in X/G$ be the projection. The triple $(X/G, d_{X/G}, p_\ast m)$ is an m.m.~space and is called the {\it orbit space} of $X$ for $G$. 

An action $G \times X \ni (g, x) \mapsto gx \in X$ is said to be {\it mm-isomorphic} if for every $g \in G$, the map $X \ni x \mapsto gx \in X$ is an isometry preserving the measure $m$. 

The family $\mathcal{F} :=\{ p^{-1}(y) \subset X \, | \, y \in X/G \}$ of orbits of an mm-isomorphic action by $G$ is a metric measure foliation on $X$.
\end{ex}

\begin{ex}[Warped product]
Warped products in the framework of m.m.~spaces are defined and studied by   \cites{Ktr, GH} et al. The warped product is an example of metric measure foliation. In the following, we define warped products along \cite{GH}. We need to assume that two m.m.~spaces defining the warped product of them are intrinsic metric spaces and at least one of them has finite measure.

Let $(Y, d_Y, m_Y), (Z, d_Z, m_Z)$ be two m.m.~spaces with intrinsic metric and $w_d, w_m : Y \to [0, +\infty)$ two bounded continuous functions on $Y$ such that $w_m \not\equiv 0$. We assume that $m_Z$ is finite. For a curve $\gamma = (\alpha, \beta)$ on $Y \times Z$ such that $\alpha, \beta$ are absolutely continuous curves on $Y, Z$ respectively, the $w_d$-length $l_w[\gamma]$ of $\gamma$ is defined by
\begin{equation}
l_w[\gamma] = \int_0^1 \sqrt{|\dot{\alpha}|(t)^2 + w_d(\alpha(t))^2|\dot{\beta}|(t)^2} \, dt.
\end{equation}
For any two points $x, x' \in Y \times Z$, we denote by $\Adm(x, x')$ the set of all curves $\gamma = (\alpha, \beta)$ on $Y \times Z$ joining $x$ and $x'$ such that $\alpha, \beta$ are absolutely continuous curve on $Y, Z$ respectively. We define a pseudo-metric $d_w$ on $Y \times Z$ by
\begin{equation}
d_w(x, x') := \inf{\{ l_w[\gamma] \mid \gamma \in \Adm(x, x')\}}
\end{equation}
for $x, x' \in Y \times Z$. The pseudo metric $d_w$ induces an equivalence relation defined by $x \sim x'$ if and only if $d_w(x, x') = 0$. The quotient space $Q := ((Y \times Z) / \sim, d_w)$ is a metric space. Since $Y, Z$ are both separable, $Q$ is also separable. We write $\hat{Q}$ as the completion of $Q$ and regard $Q$ as a subset of $\hat{Q}$. Let $q : Y \times Z \to Q$ be the quotient map. We define a Borel measure $m_w$ on $\hat{Q}$ by
\begin{equation}
m_w := q_* (w_m m_Y \otimes m_Z).
\end{equation}
Since $m_Z$ is finite and $w_m \not\equiv 0$, we see that $m_w$ is locally finite and non-trivial. We define the {\it warped product} $Y \times_w Z$ of $Y$ and $Z$ for the warping functions $w = (w_d, w_m)$ by
\begin{equation}
Y \times_w Z := (\supp{m_w}, d_w, m_w).
\end{equation}

\begin{lem}\label{warpedproduct}
Let $(Y, d_Y, m_Y), (Z, d_Z, m_Z)$ be two m.m.~spaces with intrinsic metric and $w_d, w_m : Y \to [0, +\infty)$ two bounded continuous functions on $Y$ such that $w_m \not\equiv 0$. Assume that $m_Z$ is finite. Then the projection $p : Y \times_w Z \to Y$ induces a metric measure foliation.
\end{lem}

\begin{proof}
We set $C := m_Z(Z)$ and $\Omega := \{ y \in Y \mid w_m(y) > 0 \}$. It is easy to prove that $p$ is 1-Lipschitz. Let us prove that $p$ satisfies (\ref{fc}) on $\Omega$. Let $Q$ be the quotient space and $q: Y \times Z \to Q$ the quotient map in the definition of the warped product. Since $p \circ q$ coincides with the natural projection from $Y \times Z$ to $Y$, we obtain $p_\ast m_w = Cw_m m_Y$. Moreover, we define a family $\{ \mu_y\}_{y \in Y}$ of probability measures on $Y \times_w Z$ by
\begin{equation}
\mu_y := \left\{ \begin{array}{ll} {q_y}_* (C^{-1}m_Z) &  \text{if } y \in \Omega, \\ \mu_0 & \text{otherwise} \end{array} \right.
\end{equation}
for any $y \in Y$, where $q_y := q|_{\{y\} \times Z}$ is the quotient map restricted on ${\{y\} \times Z}$ and $\mu_0$ is an arbitrary probability measure on $Y \times_w Z$. Then, since $p_* m_w(Y \setminus \Omega) = 0$, we see that $\{ \mu_y\}_{y \in Y}$ is the disintegration of $m_w$ for $p$.  We write $\nu_y$ as the measure $C^{-1}m_Z$ on $\{ y \} \times Z$. We have $\mu_y = {q_y}_* \nu_y$ for any $y \in \Omega$. Given two points $y, y' \in \Omega$, we define a Borel measurable map $\psi_{yy'} : \{ y \} \times Z \to \{ y' \} \times Z$ by
\begin{equation*}
\psi_{yy'}((y, z)) = (y', z) \in \{ y' \} \times Z
\end{equation*}
for $(y, z) \in \{ y \} \times Z$. We set $\pi_{yy'} := (q_y \times (q_{y'} \circ \psi_{yy'}))_* \nu_y$ and then have $\pi_{yy'} \in \Pi(\mu_y, \mu_{y'})$. Therefore,
\begin{align*}
W_2(\mu_y, \mu_{y'})^2 & \leq \int_{X \times X} d_w(x, x')^2 \, d\pi_{yy'}(x, x') \\
& = \int_{Z} d_w((y, z), (y', z))^2 \, d\nu_y(z) = d_Y(y, y')^2.
\end{align*}
On the other hand, by Lemma \ref{pushW_q}, we have
\begin{equation*}
W_2(\mu_y, \mu_{y'}) \geq d_Y(y, y').
\end{equation*}
These imply (\ref{fc}). The proof is completed.
\end{proof}
\end{ex}

\begin{ex}
In Section 1, we consider the sequence of the $n$-dimensional spheres $S^n(r_n)$ in $\R^{n+1}$, $n = 1, 2, \ldots$ with radii $r_n > 0$. We define a map $p_n : S^n(r_n) \to \R$ by
\begin{equation*}
p_n(x) := d_{S^n(r_n)}(x, \bar{x}_n) - \frac{\pi}{2}r_n
\end{equation*}
for $x \in S^n(r_n)$, where $\bar{x}_n$ is a fixed point in $S^n(r_n)$. For each $n$, we see that the map $p_n$ induces a metric measure foliation from the following discussion.

We define an m.m.~space $I_n$ by
\begin{equation}
I_n := \left(\left[ -\frac{\pi}{2}r_n, \frac{\pi}{2}r_n\right], | \cdot |, \mathcal{L}^1 |_{\left[ -\frac{\pi}{2}r_n, \frac{\pi}{2}r_n\right]}\right)
\end{equation}
and define two continuous maps $w_d, w_m : I_n \to [0, +\infty)$ by
\begin{equation}
w_d(t) := \cos{\frac{t}{r_n}}, \quad w_m(t) := \frac{\cos^{n-1}{\frac{t}{r_n}}}{\int_{I_n} \cos^{n-1}{\frac{t}{r_n}} \, dt}.
\end{equation}
We see that $S^n(r_n)$ is isomorphic to $I_n \times_w S^{n-1}(r_n)$ as m.m.~space and $p_n$ corresponds to the projection from $S^n(r_n)$ to $I_n$ if the fixed point $\bar{x}_n \in S^n(r_n)$ corresponds to $[(- \frac{\pi}{2}\sqrt{n}, *)] \in I_n \times_w S^{n-1}(r_n)$. By Lemma \ref{warpedproduct}, the map $p_n$ induces a metric measure foliation.
\end{ex}

\subsection{Quotient space induced by metric measure foliation}

Let $\F$ be a metric measure foliation on an m.m.~space $(X, d, m)$. We denote by $(X^*, d^*, m^*)$ the quotient m.m.~space induced by $\F$ and denote by $p : X \to X^*$ the quotient map, where the measure $m^*$ is defined by $m^* := p_* m$. 

The quotient map $p : X \to X^*$  induces a nice pullback of a probability measure on $X^*$ (which is called the lift of measure in \cite{GKMS}). 

\begin{dfn}[Pullback of measure]
Let $\nu \in \Pb(X^*)$. The {\it pullback measure} $p^\ast \nu \in \Pb(X)$ of $\nu$ by $p$ is defined by
\begin{equation}\label{pullback}
(p^\ast \nu)(A) := \int_Y \mu_y(A) \, d\nu(y)
\end{equation}
for any Borel subset $A \subset X$. 
\end{dfn}

It follows from the definition of the pullback measure $p^\ast \nu$ that for any Borel measurable function $f : X \to \mathbb{R}$,
\begin{equation}
\int_X f(x) \, d(p^\ast \nu)(x) = \int_Y \int_{p^{-1}(y)} f(x) \, d\mu_y(x)d\nu(y).
\end{equation}

\begin{rem}
For a function $f : X^* \to \mathbb{R}$, the pullback function $p^\ast f : X \to \mathbb{R}$ of $f$ by $p$ is defined by $p^\ast f := f \circ p$ naturally.
\end{rem}

\begin{prop}\label{prop_pull}
Let $\nu \in \Pb(X^*)$. We have the following (1) -- (3).
\begin{enumerate}
\item $p_\ast (p^\ast \nu) = \nu$.
\item If $\nu = \rho m^*$ for a Borel measurable function $\rho: X^* \to \mathbb{R}$, then $p^\ast \nu = (p^\ast \rho)m$.
\item If $\nu$ is absolutely continuous with respect to $p_\ast m$, then
\begin{equation*}
\Ent_{m^*}(\nu) = \Ent_{m}(p^\ast \nu).
\end{equation*}
\end{enumerate}
\end{prop}

\begin{proof}
We first prove (1). Given a Borel subset $B \subset X^*$, it holds that
\begin{align*}
p_* (p^* \nu)(B) & = (p^* \nu)(p^{-1}(B)) = \int_{X^*} \mu_y(p^{-1}(B)) \, d\nu(y) \\
 & = \int_{X^*} \mathbf{1}_B (y) \, d\nu(y) = \nu(B).
\end{align*}
This means $p_\ast (p^\ast \nu) = \nu$.

We next prove (2) and (3). We assume $\nu = \rho m^*$. For any Borel subset $A \subset X$,
\begin{align*}
p^\ast \nu(A) = & \int_{X^*} \mu_y(A) \rho(y) \, dm^*(y) \\
= & \int_{X^*} \rho(y) \int_{p^{-1}(y)} \mathbf{1}_A(x) \, d\mu_y(x)dm^*(y) \\
= & \int_{X^*} \int_{p^{-1}(y)} \rho(p(x))\mathbf{1}_A(x) \, d\mu_y(x)dm^*(y) \\
= & \int_X \mathbf{1}_A(x) \rho(p(x)) \, dm(x) = (p^\ast \rho)m (A),
\end{align*}
which implies $p^\ast \nu = (p^\ast \rho)m$. Moreover, we  have 
\begin{align*}
\Ent_{m}(p^\ast \nu) = & \int_X \rho(p(x)) \log{\rho(p(x))} \, dm(x) \\
= & \int_X \rho(y) \log{\rho(y)} \, dm^*(y) = \Ent_{m^*}(\nu).
\end{align*}
The proof is completed.
\end{proof}

\begin{lem}\label{pullW_q}
Let $\nu_0, \nu_1 \in \Pb(X^*)$ and $q \in [1, \infty)$. Then we have
\begin{equation}\label{pullW_qeq}
W_q(p^\ast \nu_0, p^\ast \nu_1) = W_q(\nu_0, \nu_1).
\end{equation}
In other words, the map
\begin{equation}
p^\ast : \Pb(X^*) \ni \nu \mapsto p^\ast \nu \in \Pb(X)
\end{equation}
is isometric with respect to $W_q$.
\end{lem}

\begin{proof}
We take any two measures $\nu_0, \nu_1 \in \Pb(X^*)$ and fix them. The inequality $W_q(\nu_0, \nu_1) \leq W_q(p^\ast \nu_0, p^\ast \nu_1)$ follows from Proposition \ref{prop_pull}(1) and Lemma \ref{pushW_q}. We prove the opposite inequality. Assume that $W_q(\nu_0, \nu_1) < +\infty$. Let $\pi \in \Pb(X^* \times X^*)$ be an optimal transport plan for $W_q(\nu_0, \nu_1)$. By Aumann's measurable choice theorem (see \cite{aumann}), there exists a family $\{ \pi_{yy'} \}_{(y, y') \in X^* \times X^*}$ of probablity measures on $X \times X$ such that the map $X^* \times X^* \ni (y, y') \mapsto \pi_{yy'}(A) \in [0, 1]$ is Borel measurable for any Borel subset $A \subset X \times X$ and $\pi_{yy'}$ is an optimal transport plan for $W_q(\mu_y, \mu_{y'})$ for $\pi$-a.e. $(y, y') \in X^* \times X^*$. (In Aumann's theorem, it is easy to check the Borel measurability of 
\begin{equation*}
\left\{ (y, y', \pi) \in (X^*)^2 \times \Pb(X^2) \midd W_q(\mu_y, \mu_{y'})^q = \int_{X \times X} d(x, x')^q \, d\pi \right\},
\end{equation*}
where $\Pb(X^2)$ has the weak topology, from (\ref{mmf_eq}).) We define a measure $\tilde{\pi} \in \Pb(X \times X)$ by 
\begin{equation}
\tilde{\pi}(A) := \int_{X^* \times X^*} \pi_{yy'}(A) \, d\pi(y, y')
\end{equation}
for any Borel subset $A \subset X \times X$. We see that $\tilde{\pi}$ is a transport plan between $p^\ast \nu_0$ and $p^\ast \nu_1$. In fact, we have
\begin{align*}
{\pr_0}_\ast \tilde{\pi}(A) = & \int_{X^* \times X^*} \pi_{yy'}({\pr_0}^{-1} (A)) \, d\pi(y, y') \\
= & \int_{X^* \times X^*} \mu_y (A) \, d\pi(y, y') 
= \int_{X^*} \mu_y (A) \, d\nu_0(y) = p^\ast \nu_0(A)
\end{align*}
for any Borel subset $A \subset X$, where $\pr_{0}$ is the projection to the first coordinate. This means ${\pr_0}_\ast \tilde{\pi} = p^\ast \nu_0$ and we obtain ${\pr_1}_\ast \tilde{\pi} = p^\ast \nu_1$ in the same way. Thus $\tilde{\pi}$ is a transport plan between $p^\ast \nu_0$ and $p^\ast \nu_1$. Then, we have
\begin{align*}
W_q(p^\ast \nu_0, p^\ast \nu_1)^q \leq & \int_{X \times X} d(x, x')^q \, d\tilde{\pi}(x, x') \\
= & \int_{X^* \times X^*} \int_{X \times X} d(x, x')^q \, d\pi_{yy'}(x, x')d\pi(y, y') \\
= & \int_{X^* \times X^*} W_q(\mu_y, \mu_{y'})^q \, d\pi(y, y') \\
= & \int_{X^* \times X^*} d^*(y, y')^q \, d\pi(y, y') = W_q(\nu_0, \nu_1)^q.
\end{align*}
By this, we obtain (\ref{pullW_qeq}). The proof is completed.
\end{proof}

In \cite{GKMS}*{Theorem 8.8}, it is shown that the strong $\CD(K, \infty)$ condition, that is, $\CD(K, \infty)$ and essentially non-branching, is preserved via the metric measure foliation with bounded leaves. The following theorem claims that the (usual) $\CD(K, \infty)$ condition is preserved even if the foliation consists of unbounded leaves. Combining the following theorem and \cite{GKMS}*{Corollary 3.5} (the foliational version of \cite{GKMS}*{Corollary 3.5} is obtained by the same proof) implies the claim of  \cite{GKMS}*{Theorem 8.8} for the strong $\CD(K, \infty)$ condition.

\begin{thm}\label{stabCD}
Let $(X, d, m)$ be an m.m.~space with a metric measure foliation and let $K \in \mathbb{R}$. Assume that $X$ satisfies $\CD(K, \infty)$ and $X^*$ satisfies $\VG$. Then the quotient space $X^*$ satisfies $\CD(K, \infty)$.
\end{thm}

Note that if $X^*$ satisfies $\CD(K, \infty)$, by Proposition \ref{OntheVG}, then $X^*$ satisfies $\VG$, so that the condition $\VG$ is a natural assumption.

\begin{proof}
The outline of the proof is the same as that of \cite{GKMS}*{Theorem 8.8}.

Since $X^*$ satisfies $\VG$, there exist $\bar{y} \in X^*$ and $C > 0$ such that
\begin{equation}
\int_{X^*} e^{-C^2 d^*(y, \bar{y})^2} \, dm^*(y) < + \infty.
\end{equation}
We define $V : X \to [0, \infty)$ by $V(x) := C d^*(p(x), \bar{y})$ for $x \in X$.

We take any two measures $\nu_0, \nu_1 \in \Pb_2(X^*) \cap D(\Ent_{m^*})$. Let $\rho_0, \rho_1$ be the densities of $\nu_0, \nu_1$ with respect to $m^*$ respectively, that is, $\nu_i = \rho_i m^*$ for $i = 0, 1$. Then, by Proposition \ref{prop_pull}, $p^\ast \nu_0, p^\ast \nu_1$ both belong to $\Pb_V(X) \cap D(\Ent_m)$. Moreover, by Lemma \ref{pullW_q}, we have
\begin{equation*}
W_2(p^\ast \nu_0, p^\ast \nu_1) = W_2(\nu_0, \nu_1) < + \infty.
\end{equation*}
Thus, by Lemma \ref{CDlem}, there exists a $W_2$-geodesic $\mu : [0,1] \ni t \mapsto \mu_t \in \Pb_V(X)$ joining $p^\ast \nu_0$ and $p^\ast \nu_1$ such that for any $t \in [0, 1]$,
\begin{equation}\label{pullKconvex}
\Ent_m (\mu_t) \leq (1-t) \Ent_m (p^\ast \nu_0) + t \Ent_m (p^\ast \nu_1) - \frac{K}{2} t(1-t) W_2(p^\ast \nu_0, p^\ast \nu_1)^2.
\end{equation}
We set $\nu_t := p_\ast \mu_t \in \Pb_2(X^*)$ for any $t \in (0, 1)$. Let us prove that $\nu : [0,1] \ni t \mapsto \nu_t \in \Pb_2(X^*)$ is a $W_2$-geodesic joining $\nu_0$ and $\nu_1$ satisfying (\ref{Kconvex}). By Proposition \ref{prop_pull} (1), we have
\begin{equation*}
p_\ast \mu_i = p_\ast (p^\ast \nu_i) = \nu_i
\end{equation*}
for $i = 0, 1$. Combining Lemma \ref{pushW_q}, the definition of $W_2$-geodesic, and Lemma \ref{pullW_q} yields that for any $s, t \in [0, 1]$,
\begin{equation*}
W_2(\nu_s, \nu_t) \leq W_2(\mu_s, \mu_t) = |s - t| W_2(\mu_0, \mu_1) = |s - t| W_2(\nu_0, \nu_1).
\end{equation*}
On the other hand, by the triangle inequality,
\begin{align*}
W_2(\nu_s, \nu_t) & \geq W_2(\nu_0, \nu_1) - W_2(\nu_0, \nu_s) - W_2(\nu_t, \nu_1) \\
& \geq W_2(\nu_0, \nu_1) - s W_2(\nu_0, \nu_1) - (1 - t) W_2(\nu_0, \nu_1) \\
& = (t - s) W_2(\nu_0, \nu_1),
\end{align*}
which implies that for any $s, t \in [0, 1]$, 
\begin{equation*}
W_2(\nu_s, \nu_t) = |s - t| W_2(\nu_0, \nu_1).
\end{equation*}
Thus, $\nu : [0, 1] \ni t \mapsto \nu_t \in \Pb_2(X^*)$ is a $W_2$-geodesic joining $\nu_0$ and $\nu_1$. Moreover, combining (\ref{pullKconvex}), Proposition \ref{prop_pull}, and Lemma \ref{pullW_q} yields that for any $t \in [0, 1]$,
\begin{align*}
& \Ent_{m^*}(\nu_t) \leq \Ent_m(\mu_t) \\
& \leq (1-t) \Ent_m (p^\ast \nu_0) + t \Ent_m (p^\ast \nu_1) - \frac{K}{2} t(1-t) W_2(p^\ast \nu_0, p^\ast \nu_1)^2 \\
& = (1-t) \Ent_{m^*} (\nu_0) + t \Ent_{m^*} (\nu_1) - \frac{K}{2} t(1-t) W_2(\nu_0, \nu_1)^2,
\end{align*}
which implies that $\nu$ satisfies (\ref{Kconvex}). Therefore, we see that $X^*$ satisfies $\CD(K, \infty)$. The proof is completed.
\end{proof}

\begin{thm}\label{stabCh}
Let $(X, d, m)$ be an m.m.~space with a metric measure foliation and let $q \in (1, \infty)$. Then we have
\begin{equation}
\Ch^{X^*}_q(f) = \Ch^X_q(p^* f)
\end{equation}
for any $f \in L^q(X^*, m^*)$.
\end{thm}

This theorem corresponds to \cite{GKMS}*{Proposition 8.9 (4)}. Owing to an approximation by Lipschitz functions, some assumptions for the Sobolev space is required in \cite{GKMS}*{Proposition 8.9}. Using a new technique, we remove the assumptions. We prove the following corollary using Theorem \ref{stabCh} before we prove Theorem \ref{stabCh}.

\begin{cor}\label{stabHilb}
Let $(X, d, m)$ be an m.m.~space with a metric measure foliation. Assume that $X$ is infinitesimally Hilbertian. Then the quotient space $X^*$ is infinitesimally Hilbertian.
\end{cor}

\begin{proof}
We take any two functions $f, g \in L^2(X^*, m^*)$. Let us prove that 
\begin{equation*}
\Ch_2^{X^*}(f + g) + \Ch_2^{X^*}(f - g) = 2 \Ch_2^{X^*}(f) + 2 \Ch_2^{X^*}(g).
\end{equation*}
Since their pullback functions $p^* f, p^* g$ belong to $L^2(X, m_X)$ and the functional $\Ch_2^X$ is quadratic, it holds that
\begin{equation*}
\Ch_2^X(p^* f + p^* g) + \Ch_2^X(p^* f - p^* g) = 2 \Ch_2^X(p^* f) + 2 \Ch_2^X(p^* g).
\end{equation*}
Thus, by Theorem \ref{stabCh}, we have
\begin{align*}
& \Ch_2^{X^*}(f + g) + \Ch_2^{X^*}(f - g) = \Ch_2^X(p^* f + p^* g) + \Ch_2^X(p^* f - p^* g) \\
& = 2 \Ch_2^X(p^* f) + 2 \Ch_2^X(p^* g) = 2 \Ch_2^{X^*}(f) + 2 \Ch_2^{X^*}(g).
\end{align*}
This completes the proof.
\end{proof}

Combining Theorem \ref{stabCD} and Corollay \ref{stabHilb} proves the following.

\begin{cor}\label{stabRCD}
Let $(X, d, m)$ be an m.m.~space with a metric measure foliation and let $K \in \mathbb{R}$. Assume that $X$ satisfies $\RCD(K, \infty)$ and $X^*$ satisfies $\VG$. Then the quotient space $X^*$ satisfies $\RCD(K, \infty)$.
\end{cor}

We prove Theorem \ref{stabCh}. Let $(X, d, m)$ be an m.m.~space with a metric measure foliation and let $q \in (1, \infty)$. We write $\widetilde{X}$, $\widetilde{X^*}$ as $C([0, 1]; X)$, $C([0, 1]; X^*)$ respectively for simplicity. Moreover, we set a map $\tilde{p}: \widetilde{X} \to \widetilde{X^*}$ by $\tilde{p}(\xi) := p \circ \xi$ for any curve $\xi \in \widetilde{X}$. The map $\tilde{p}$ is 1-Lipschitz with respect to the uniform distance. We first obtain the following proposition.

\begin{prop}\label{pushwug}
Let $(X, d, m)$ be an m.m.~space with a metric measure foliation and $q \in (1, \infty)$. Then, for any $f \in W^{1,q}(X^*, d^*, m^*)$, the pullback $p^* f$ belongs to $W^{1,q}(X, d, m)$ and
\begin{equation}\label{pushwug_eq}
|D(p^* f)|_w(x) \leq |Df|_w(p(x))
\end{equation}
holds for $m$-a.e.~$x \in X$. In particular, for any $f \in L^q(X^*, m^*)$, we have 
\begin{equation}\label{push_Ch}
\Ch_q^X(p^* f) \leq \Ch_q^{X^*}(f).
\end{equation}
\end{prop}

\begin{rem}
Proposition \ref{pushwug} does not need the foliational structure. We obtain the same result for a 1-Lipschitz map $p : X \to Y$ between two m.m.~spaces $X$ and $Y$ satisfying $p_* m_X = m_Y$. 
\end{rem}

\begin{proof}
We take any $f \in W^{1,q}(X^*, d^*, m^*)$ and fix it. For the proofs of $p^* f \in W^{1,q}(X, d, m)$ and (\ref{pushwug_eq}), it suffices to prove
\begin{align*}
& \int_{\widetilde{X}} |(p^* f)(\xi(1)) - (p^* f)(\xi(0))| \, d\pi(\xi) \\ 
& \leq \int_{\widetilde{X}}\int_0^1 |Df|_w(p(\xi(t))) |\dot{\xi}|(t) \, dtd\pi(\xi)
\end{align*}
for any $q^*$-test plan $\pi \in \Pb(\widetilde{X})$ on $X$, where $q^*$ is the conjugate exponent of $q$. We take any $q^*$-test plan $\pi \in \Pb(\widetilde{X})$ and fix it. Then the measure $\tilde{p}_* \pi \in \Pb(\widetilde{X^*})$ is a $q^*$-test plan on $X^*$. In fact, since $\pi$ is a $q^*$-test plan, there exists a constant $C > 0$ such that ${e_t}_\ast \pi \leq Cm$ for any $t \in [0, 1]$. Since $e_t \circ \tilde{p} = p \circ e_t$ for each $t \in [0, 1]$, we have
\begin{equation*}
{e_t}_* \tilde{p}_* \pi = p_* {e_t}_* \pi \leq p_* (Cm) = Cm^*.
\end{equation*}
Moreover, it holds that
\begin{equation*}
\int_{\widetilde{X^*}} \mathcal{E}_{q^*}[\gamma] \, d(\tilde{p}_* \pi)(\gamma) = \int_{\widetilde{X}} \mathcal{E}_{q^*}[p \circ \xi] \, d\pi(\xi) \leq \int_{\widetilde{X}} \mathcal{E}_{q^*}[\xi] \, d\pi(\xi) < +\infty,
\end{equation*}
which implies that $\tilde{p}_* \pi$ is a $q^*$-test plan on $Y$. Thus, by the definition of $|Df|_w$, it holds that
\begin{align*}
& \int_{\widetilde{X^*}} |f(\gamma(1)) - f(\gamma(0))| \, d(\tilde{p}_* \pi)(\gamma) \\
& \leq \int_{\widetilde{X^*}} \int_0^1 |Df|_w(\gamma(t)) |\dot{\gamma}|(t) \, dt d(\tilde{p}_* \pi)(\gamma).
\end{align*}
This implies 
\begin{align*}
& \int_{\widetilde{X}} |(p^\ast f)(\xi(1)) - (p^* f)(\xi(0))| \, d\pi(\xi) \\
& = \int_{\widetilde{X^*}} |f(\gamma(1)) - f(\gamma(0))| \, d(\tilde{p}_* \pi)(\gamma) \\
& \leq \int_{\widetilde{X^*}} \int_0^1 |Df|_w(\gamma(t)) |\dot{\gamma}|(t) \, dt d(\tilde{p}_* \pi)(\gamma) \\
& = \int_{\widetilde{X}} \int_0^1 |Df|_w(p(\xi(t))) |\dot{(p \circ \xi)}|(t) \, dt d\pi(\xi) \\
& \leq \int_{\widetilde{X}} \int_0^1 |Df|_w(p(\xi(t))) |\dot{\xi}|(t) \, dt d\pi(\xi).
\end{align*}
The proof is completed.
\end{proof}

It is sufficient to show the opposite inequality of (\ref{push_Ch}) for the proof of Theorem \ref{stabCh}. We prepare the following to prove it. Let $\{\mu_y\}_{y \in X^*}$ be canonical disintegration of $m$ for the quotient map $p$. We fix $q \in (1, \infty)$ and denote $q^*$ the conjugate exponent of $q$.

We take any curve $\gamma \in AC^{q^*}([0, 1]; X^*) \subset \widetilde{X^*}$ and consider the curve $\mu_\gamma : [0, 1] \ni t \mapsto \mu_{\gamma(t)} \in \Pb(X)$. Since we have
\begin{equation*}
W_{q^*}(\mu_{\gamma(s)}, \mu_{\gamma(t)}) = d_Y(\gamma(s), \gamma(t))
\end{equation*}
for any $s, t \in [0, 1]$, the curve $\mu_\gamma$ belongs to $AC^{q^*}([0, 1]; (\Pb(X), W_{q^*}))$ and $|\dot{\mu_\gamma}|(t) = |\dot{\gamma}|(t)$ holds for $\mathcal{L}^1$-a.e.~$t \in [0, 1]$. Thus, by Proposition \ref{lisinithm}, there exists a measure $\eta_\gamma \in \Pb(\widetilde{X})$ such that 
\begin{align}
& {e_t}_\ast \eta_\gamma = \mu_{\gamma(t)}  \text{ for any } t \in [0, 1], \label{lift_disint} \\
& \int_{\widetilde{X}} \mathcal{E}_{q^*} [\xi] \, d\eta_\gamma(\xi) = \int_0^1 |\dot{\mu_\gamma}|(t)^{q^*} \, dt = \int_0^1 |\dot{\gamma}|(t)^{q^*} \, dt = \mathcal{E}_{q^*} [\gamma]. \label{lift_energy}
\end{align}

Since $e_t \circ \tilde{p} = p \circ e_t$ for each $t \in [0, 1]$, we have $(e_t)_* \tilde{p}_* \eta_\gamma = \delta_{\gamma(t)}$ for each $t \in [0, 1]$. This implies  $\tilde{p}_* \eta_\gamma = \delta_\gamma \in \Pb(\widetilde{X^*})$. By Aumann's measurable choice theorem (see \cite{aumann}), whenever we take a probability measure $\pi$ on $\widetilde{X^*}$, there exists a family $\{\eta_{\gamma}\}_{\gamma \in \widetilde{X^*}}$ of probability measures on $\widetilde{X}$ such that the map $\widetilde{X^*} \ni \gamma \mapsto \eta_{\gamma}(A) \in [0, 1]$ is Borel measurable for any Borel subset $A \subset \widetilde{X}$ and, for $\pi$-a.e. $\gamma \in \widetilde{X^*}$, the measure $\eta_{\gamma}$ satisfies (\ref{lift_disint}) and (\ref{lift_energy}).

From this discussion, we regard the family $\{\eta_\gamma\}_{\gamma \in \widetilde{X^*}} \subset \Pb(\widetilde{X})$ as a ``nice'' lift of the disintegration $\{\mu_y\}_{y \in X^*} \subset \Pb(X)$ of $m$. Using this family $\{\eta_\gamma\}_{\gamma \in \widetilde{X^*}}$, we are able to lift all test plans on $\widetilde{X^*}$ to $\widetilde{X}$.

\begin{prop}\label{pull_testplan}
Let $\pi \in \Pb(\widetilde{X^*})$ be a $q^*$-test plan on $X^*$. We define a measure $\tilde{p}^* \pi \in \Pb(\widetilde{X})$ by 
\begin{equation}\label{pull_testplan_eq}
\tilde{p}^* \pi(A) := \int_{\widetilde{X^*}} \eta_\gamma(A) \, d\pi(\gamma) 
\end{equation}
for any Borel subset $A \subset \widetilde{X}$. Then $\tilde{p}^* \pi$ is a $q^*$-test plan on $X$ and satisfies $\tilde{p}_\ast (\tilde{p}^* \pi) = \pi$.
\end{prop}

\begin{proof}
Let $\pi \in \Pb(\widetilde{X^*})$ be a $q^*$-test plan on $X^*$. By the definition of $q^*$-test plan, there exists a constant $C > 0$ such that ${e_t}_\ast \pi \leq Cm^*$ for any $t \in [0, 1]$. We first prove that
\begin{equation}\label{finite_energy}
\int_{\widetilde{X}} \mathcal{E}_{q^*}[\xi] \, d(\tilde{p}^* \pi)(\xi) < +\infty.
\end{equation}
By the definition of $\tilde{p}^* \pi$ and (\ref{lift_energy}), we see that
\begin{equation*}
\int_{\widetilde{X}} \mathcal{E}_{q^*}[\xi] \, d(\tilde{p}^* \pi)(\xi) = \int_{\widetilde{X^*}} \int_{\widetilde{X}} \mathcal{E}_{q^*}[\xi] \, d\eta_\gamma(\xi)d\pi(\gamma) = \int_{\widetilde{X^*}} \mathcal{E}_{q^*}[\gamma] \, d\pi(\gamma).
\end{equation*}
We obtain (\ref{finite_energy}). We next prove that ${e_t}_\ast (\tilde{p}^* \pi) \leq Cm$ for any $t \in [0, 1]$. We take any $t \in [0, 1]$. Then, for any Borel subset $A \subset X$,
\begin{align*}
{e_t}_\ast (\tilde{p}^* \pi)(A) & = \int_{\widetilde{X^*}} \eta_\gamma(e_t^{-1}(A)) \, d\pi(\gamma) = \int_{\widetilde{X^*}} \mu_{\gamma(t)}(A) \, d\pi(\gamma) \\
& = \int_{X^*} \mu_{y}(A) \, d({e_t}_* \pi)(y) \leq C \int_{X^*} \mu_{y}(A) \, dm^*(y) = Cm(A),
\end{align*}
which implies ${e_t}_\ast (\tilde{p}^* \pi) \leq Cm$. Thus, we see that $\tilde{p}^* \pi$ is a $q^*$-test plan on $X$. Moreover, we see that $\tilde{p}^* \pi$ satisfies $\tilde{p}_\ast (\tilde{p}^* \pi) = \pi$ in the same way as in the proof of Proposition \ref{prop_pull} (1). The proof is completed.
\end{proof}

\begin{lem}\label{compare_wug}
Let $f \in W^{1,q}(X, d, m)$. We define $g \in L^q(X^*, m^*)$ by
\begin{equation}\label{push_fct}
g(y) := \int_{X} f(x) \, d\mu_y(x)
\end{equation}
for $y \in X^*$. Then $g$ belongs to $W^{1,q}(X^*, d^*, m^*)$ and
\begin{equation}\label{compare_wug_eq}
|Dg|_w(y)^q \leq \int_X |Df|_w(x)^q \, d\mu_y(x)
\end{equation}
holds for $m^*$-a.e.~$y \in X^*$.
\end{lem}

\begin{rem}
Given a function $h \in W^{1,q}(X^*, d^*, m^*)$, we apply Lemma \ref{compare_wug} to the function $f := p^* h$. Then $g$ coincides with $h$ and
\begin{equation}\label{compare_wug_eq_pull}
|Dh|_w(y)^q \leq \int_X |D(p^* h)|_w(x)^q \, d\mu_y(x)
\end{equation}
holds for $m^*$-a.e.~$y \in X^*$. For the proof of Theorem \ref{stabCh}, it suffices to prove (\ref{compare_wug_eq_pull}). We obtain the stronger formula (\ref{compare_wug_eq}).
\end{rem}

\begin{proof}[{\bf Proof of Lemma \ref{compare_wug}}]
We take any $f \in W^{1,q}(X, d, m)$ and define $g \in L^q(X^*, m^*)$ by (\ref{push_fct}). For the proof of $g \in W^{1,q}(X^*, d^*, m^*)$ and (\ref{compare_wug_eq}), it is sufficient to prove 
\begin{equation}\label{pull_wug_def_eq}
\begin{split}
& \int_{\widetilde{X^*}} |g(\gamma(1)) - g(\gamma(0))| \, d\pi(\gamma) \\
& \leq \int_{\widetilde{X^*}} \int_0^1 \left(\int_X |Df|_w(x)^q \, d\mu_{\gamma(t)}(x) \right)^{\frac{1}{q}} |\dot{\gamma}|(t) \, dt d\pi(\gamma)
\end{split}
\end{equation}
for any $q^*$-test plan $\pi \in \Pb(\widetilde{X^*})$ on $X^*$. We take any $q^*$-test plan $\pi \in \Pb(\widetilde{X^*})$ on $X^*$ and fix it. By Proposition \ref{pull_testplan}, the measure $\tilde{p}^* \pi$ defined by (\ref{pull_testplan_eq}) is a $q^*$-test plan on $X$. By the definition of $|Df|_w$, we have
\begin{equation*}
\int_{\widetilde{X}} |f(\xi(1)) - f(\xi(0))| \, d(\tilde{p}^* \pi)(\xi) \leq \int_{\widetilde{X}} \int_0^1 |Df|_w(\xi(t)) |\dot{\xi}|(t) \, dt d(\tilde{p}^*\pi)(\xi).
\end{equation*}
Moreover, taking ${e_t}_* \eta_\gamma = \mu_{\gamma(t)}$ into account, we have
\begin{equation*}
g(\gamma(t)) = \int_X f(x) \, d\mu_{\gamma(t)}(x) = \int_{\widetilde{X}} f(\xi(t)) \, d\eta_\gamma(\xi).
\end{equation*}
Therefore, 
\begin{align*}
& \int_{\widetilde{X^*}} |g(\gamma(1)) - g(\gamma(0))| \, d\pi(\gamma) \\
& \leq \int_{\widetilde{X^*}} \int_{\widetilde{X}} |f(\xi(1)) - f(\xi(0))| \, d\eta_\gamma(\xi)d\pi(\gamma) \\
& = \int_{\widetilde{X}} |f(\xi(1)) - f(\xi(0))| \, d(\tilde{p}^* \pi)(\xi) \\
& \leq \int_{\widetilde{X}} \int_0^1 |Df|_w(\xi(t)) |\dot{\xi}|(t) \, dt d(\tilde{p}^*\pi)(\xi) \\
& = \int_{\widetilde{X^*}} \int_{\widetilde{X}} \int_0^1 |Df|_w(\xi(t)) |\dot{\xi}|(t) \, dt d\eta_\gamma(\xi)d\pi(\gamma) \\
& \leq \int_{\widetilde{X^*}} \int_0^1 \left( \int_{\widetilde{X}}  |Df|_w(\xi(t))^q \, d\eta_\gamma(\xi) \right)^\frac{1}{q} \left( \int_{\widetilde{X}} |\dot{\xi}|(t)^{q^*} \, d\eta_\gamma(\xi) \right)^\frac{1}{q^*} \, dt d\pi(\gamma) \\
& = \int_{\widetilde{X^*}} \int_0^1 \left( \int_X  |Df|_w(x)^q \, d\mu_{\gamma(t)}(x) \right)^\frac{1}{q} \left( \int_{\widetilde{X}} |\dot{\xi}|(t)^{q^*} \, d\eta_\gamma(\xi) \right)^\frac{1}{q^*} \, dt d\pi(\gamma).
\end{align*}
In order to prove (\ref{pull_wug_def_eq}), it suffices to prove  
\begin{equation}\label{ptwise_energy}
\int_{\widetilde{X}} |\dot{\xi}|(t)^{q^*} \, d\eta_\gamma(\xi) = |\dot{\gamma}|(t)^{q^*}
\end{equation}
for $\pi$-a.e. $\gamma \in \widetilde{X^*}$ and $\mathcal{L}^1$-a.e.~$t \in [0, 1]$. Since $p$ is 1-Lipschitz, we have
\begin{equation*}
\int_{\widetilde{X}} |\dot{\xi}|(t)^{q^*} \, d\eta_\gamma(\xi) \geq \int_{\widetilde{X}} |\dot{(p \circ \xi)}|(t)^{q^*} \, d\eta_\gamma(\xi) = |\dot{\gamma}|(t)^{q^*}
\end{equation*}
for $\pi$-a.e. $\gamma \in \widetilde{X^*}$ and $\mathcal{L}^1$-a.e.~$t \in [0, 1]$. On the other hand, by (\ref{lift_energy}),
\begin{equation*}
\int_0^1 \left( \int_{\widetilde{X}} |\dot{\xi}|(t)^{q^*} \, d\eta_\gamma(\xi) - |\dot{\gamma}|(t)^{q^*} \right) \, dt = 0
\end{equation*}
holds for $\pi$-a.e.~$\gamma \in \widetilde{X^*}$. These imply (\ref{ptwise_energy}). Thus we obtain (\ref{pull_wug_def_eq}). This completes the proof.
\end{proof}

\begin{proof}[{\bf Proof of Theorem \ref{stabCh}}]
We take any $f \in L^q(X^*, m^*)$. By Proposition \ref{pushwug}, we have
\begin{equation*}
\Ch_q^X(p^* f) \leq \Ch_q^{X^*}(f).
\end{equation*}
If $\Ch_q^X(p^* f) = +\infty$, then we obtain $\Ch_q^X(p^* f) = \Ch_q^{X^*}(f)$. We consider the case of $p^* f \in W^{1,q}(X, d, m)$. Then, by Lemma \ref{compare_wug},
\begin{equation*}
|Df|_w(y)^q \leq \int_X |D(p^* f)|_w(x)^q \, d\mu_y(x)
\end{equation*}
holds for $m^*$-a.e. $y \in X^*$. Therefore,
\begin{align*}
\Ch_q^{X^*}(f) & = \frac{1}{q} \int_{X^*} |Df|_w(y)^q \, dm^*(y) \\
& \leq \frac{1}{q} \int_{X^*} \int_X |D(p^* f)|_w(x)^q \, d\mu_y(x)dm^*(y) \\
& = \frac{1}{q} \int_X |D(p^* f)|_w(x)^q \, dm(x) = \Ch_q^X(p^* f),
\end{align*}
which implies $\Ch_q^X(p^* f) = \Ch_q^{X^*}(f)$. The proof is completed.
\end{proof}

%% file: sec4_convergence.tex
\section{Convergence under the metric measure foliation}
In this section, we study the stability of the curvature-dimension condition and the convergences of the $q$-Cheeger energy functionals $\Ch_q^{X_n}$ and the descending slopes $\DEnt{m_n}$ (Definition \ref{DEnt}) of the relative entropies $\Ent_{m_n}$ for a sequence of p.m.m.~spaces $(X_n, d_n, m_n, \bar{x}_n)$ with a metric measure foliation. 

\subsection{Variational convergence of $q$-Cheeger energy functionals}

Our goal in this subsection is to prove Theorem \ref{intro_main}.  We denote by $\N$ the set of positive integers.

\begin{dfn}\label{pmGconv}
Let $\{ (X_n, d_n, m_n, \bar{x}_n)\}_{n \in \mathbb{N}}$ be a sequence of p.m.m. spaces and $(Y, d, m, \bar{y})$ a p.m.m.~space. We say that $\{X_n\}_{n \in \mathbb{N}}$ pmG-{\it converges} to $Y$ if there exist a complete separable metric space $(Z, d_Z)$ and isometric embeddings $\iota_n: X_n \to Z$, $\iota: Y \to Z$ such that ${\iota_n}_* m_n$ converges weakly to $\iota_* m$ in $\Mloc(Z)$ and $\iota_n(\bar{x}_n)$ converges to $\iota(\bar{y})$ in $Z$ as $n \to \infty$.
\end{dfn}

The following result was obtained by Gigli-Mondino-Savar\'e \cite{GMS}.

\begin{thm}[\cite{GMS}*{Theorem 4.9 and 7.2}]\label{GMSCD}
Let $\{(X_n, d_n, m_n, \bar{x}_n)\}_{n \in \mathbb{N}}$ be a sequence of p.m.m.~spaces and $(Y, d, m, \bar{y})$ a p.m.m.~space. Assume that $X_n$ pmG-converges to $Y$ and each $X_n$ satisfies $\CD(K, \infty)$ ($\resp{\RCD(K, \infty)}$) for a real number $K \in \mathbb{R}$. Then, $Y$ also satisfies $\CD(K, \infty)$ ($\resp{\RCD(K, \infty)}$).
\end{thm}

Combining this theorem and Theorem \ref{stabCD} and Corollary \ref{stabRCD} proves Theorem \ref{intro_main} (1) and (3) directly.

\begin{proof}[{\bf Proof of Theorem \ref{intro_main} (1) and (3)}]
We only prove (1). By Theorem \ref{stabCD}, if each $X_n$ satisfies $\CD(K, \infty)$, then $X_n^*$ satisfies $\CD(K, \infty)$. Since $X_n^*$ pmG-converges to $Y$ and by Theorem \ref{GMSCD}, the space $Y$ also satisfies $\CD(K, \infty)$. The proof is completed. We obtain (3) in the same way using Corollary \ref{stabRCD}.
\end{proof}

In order to prove Theorem \ref{intro_main} (2) and (4), we need to define the appropriate $\Gamma$-convergence and Mosco convergence of Cheeger energy functionals in our setting.

Let $\{(X_n, d_n, m_n, \bar{x}_n)\}_{n \in \mathbb{N}}$ be a sequence of p.m.m.~spaces and $(Y, d, m, \bar{y})$ a p.m.m.~space. From now on, we assume that each $X_n$ has a metric measure foliation and its quotient space $X_n^*$ pmG-converges to $Y$ as $n \to \infty$. Let $Z$ be a complete separable metric space in Definition \ref{pmGconv} associated with $X_n^* \pmG Y$. We are able to regard $X_n^*$, $Y$ as the subsets on $Z$ via the isometric embeddings. Furthermore, we denote by $p_n : X_n \to X_n^*$ the quotient map and assume that every p.m.m.~space $X_n$ satisfies $\VG$.

\begin{dfn}\label{Lweakconv}
Let $\{q_n\}_{n \in \N} \subset (1, \infty)$ be a sequence of real numbers convergent to a real number $q \in (1, \infty)$. Let $f_n \in L^{q_n}(X_n, m_n)$ for each $n \in \mathbb{N}$ and $f \in L^q(Y, m)$. We say that $f_n$ {\it $L^{q_n}$-weakly converges to} $f$ if for any function $\varphi \in \Cbs(Z)$, 
\begin{equation}\label{Lweak}
\lim_{n \to \infty} \int_{X_n} \varphi(p_n(x)) f_n(x) \, dm_n(x) = \int_{Y} \varphi(y) f(y) \, dm(y)
\end{equation}
holds, and 
\begin{equation}\label{Lbdd}
\limsup_{n \to \infty} \| f_n \|_{L^{q_n}(X_n, m_n)} < +\infty
\end{equation}
holds. 
\end{dfn}

\begin{dfn}\label{Lstrconv}
Let $\{q_n\}_{n \in \N} \subset [1, \infty)$ be a sequence of real numbers convergent to a real number $q \in [1, \infty)$ as $n \to \infty$. Let $f_n \in L^{q_n}(X_n, m_n)$ for each $n \in \mathbb{N}$ and $f \in L^q(Y, m)$.
\begin{enumerate}
\item In the case that $q_n = q > 1$, we say that $f_n$ {\it $L^q$-strongly converges to} $f$ if $f_n$ $L^q$-weakly converges to $f$ and it holds that
\begin{equation}\label{normconv}
\lim_{n \to \infty} \| f_n \|_{L^q(X_n, m_n)} = \| f \|_{L^{q}(Y, m)}.
\end{equation}
\item In the case that $q_n = q =1$, we say that $f_n$ {\it $L^1$-strongly converges to} $f$ if $\sigma \circ f_n$ $L^2$-strongly converges to $\sigma \circ f$, where $\sigma(t) := \mathrm{sign}(t)\sqrt{|t|}$ is the signed root.
\item In the case that $q_n \to q  > 1$, we say that $f_n$ {\it $L^{q_n}$-strongly converges to} $f$ if $f_n$ $L^{q_n}$-weakly converges to $f$ and, for any $\varepsilon > 0$, there exists a decomposition $f_n = \alpha_n + \beta_n$ such that
\begin{itemize}
\item $\sup_{n} \| \alpha_n \|_{L^\infty(X_n, m_n)} < +\infty$ and $\alpha_n$ $L^1$-strongly convergent,
\item $\sup_{n} \| \beta_n \|_{L^{q_n}(X_n, m_n)} < \varepsilon$.
\end{itemize}
\end{enumerate}
\end{dfn}

\begin{rem}
In the case that $(X_n, d_n, m_n, \bar{x}_n) = (Y, d, m, \bar{y})$ and $p_n = {\rm id}_Y$ for each $n \in \mathbb{N}$, the $L^q$-convergence in Definition \ref{Lweakconv} and \ref{Lstrconv} is equivalent to the usual $L^q$-convergence on $L^q(Y, m)$. In the  case that $X_n$ pmG-conveges to $Y$, Definition \ref{Lweakconv} and \ref{Lstrconv} coincide with \cite{GMS}*{Section 6.1} and \cite{AH}*{Section 3}.
\end{rem}

A discussion in the same way as \cite{GMS}, \cite{gradflow}*{Section 5.4}, \cite{AH}*{Section 3} yields that the $L^{q_n}$-convergence in Definition \ref{Lweakconv} and \ref{Lstrconv} has some basic properties as follows. From now on, we assume that $\{q_n\}_{n \in \N} \subset (1, \infty)$ converges to $q \in (1, \infty)$.

\begin{prop}\label{weaklscnorm}
Let $f_n \in L^{q_n}(X_n, m_n), f \in L^q(Y, m)$ and assume that $f_n$ $L^{q_n}$-weakly converges to $f$. Then, it holds that
\begin{equation}\label{lscnorm}
\liminf_{n \to \infty} \| f_n \|_{L^{q_n}(X_n, m_n)} \geq \| f \|_{L^{q}(Y, m)}.
\end{equation}
Moreover, any sequence $f_n \in L^{q_n}(X_n, m_n)$ with (\ref{Lbdd}) has a $L^{q_n}$-weakly converging subsequence.
\end{prop}

This proposition is proved by the way in \cite{gradflow}*{Section 5.4} directly. On the other hand, there is the following easier proof by admitting this property of the $L^{q_n}$-convergence in the pmG-convergent case in \cite{AH}*{Proposition 3.1}.

\begin{proof}[{\bf Proof of Proposition \ref{weaklscnorm}}]
We first prove (\ref{lscnorm}). We take any functions $f_n \in L^{q_n}(X_n, m_n)$, $f \in L^q(Y, m)$ and assume that $f_n$ $L^{q_n}$-weakly converges to $f$. For each $n$, we define a function $g_n \in L^{q_n}(X_n^*, m_n^*)$  by
\begin{equation*}
g_n(z) := \int_{X_n} f_n(x) \, d\mu_z^n(x),
\end{equation*}
where $\{ \mu_z^n \}_{z \in X_n^*}$ is the disintegration of $m_n$ for $p_n$. Actually, we have
\begin{equation} \label{avr_norm}
\begin{split}
\int_{X_n^*} g_n(z)^{q_n} \, dm_n^*(z) = & \int_{X_n^*} \left( \int_{X_n} f_n(x) \, d\mu_z^n (x) \right)^{q_n} dm_n^*(z) \\
\leq & \int_{X_n} f_n(x)^{q_n} \, dm_n(x),
\end{split}
\end{equation}
which implies $g_n \in L^{q_n}(X_n^*, m_n^*)$. For any $\varphi \in \Cbs(Z)$, we have
\begin{align*}
\int_{X_n^*} \varphi(z) g_n(z) \, dm_n^*(z) = & \int_{X_n^*} \int_{X_n} \varphi(z) f_n(x) \, d\mu_z^n (x) dm_n^*(z) \\
= & \int_{X_n} \varphi(p_n(x)) f_n(x) \, dm_n(x).
\end{align*}
Thus, since $f_n$ $L^{q_n}$-weakly converges to $f$ as $n \to \infty$, we obtain 
\begin{equation*}
\lim_{n \to \infty} \int_{X_n^*} \varphi(z) g_n(z) \, dm_n^*(z) = \int_Y \varphi(y) f(y) \, dm(y),
\end{equation*}
which implies that $g_n$ $L^{q_n}$-weakly converges to $f$ on $Z$ if we regard $g_n$ and $f$ as the functions on $Z$. By the lower semicontinuity of $L^{q_n}$-norm in the pmG-convergent case, we obtain
\begin{equation} \label{pmg_lsc}
\liminf_{n \to \infty} \| g_n \|_{L^{q_n}(X_n^*, m_n^*)} \geq \| f \|_{L^q(Y, m)}.
\end{equation}
Combining (\ref{avr_norm}) and (\ref{pmg_lsc}) leads to (\ref{lscnorm}). 

We next prove the weak compactness of the $L^{q_n}$-bounded sequence. We take any functions $f_n \in L^{q_n}(X_n, m_n)$ satisfying (\ref{Lbdd}). We define $g_n \in L^{q_n}(X_n^*, m_n^*)$ as in the same way as above. By (\ref{avr_norm}), the sequence $\{ g_n \}_{n \in \N}$ also satisfies (\ref{Lbdd}). By the weak compactness in the pmG-convergent case, there exists $f \in L^q(Y, m)$ such that $g_n$ $L^{q_n}$-weakly converges to $f$ on $Z$ as $n \to \infty$. In the same way as above, we see that $f_n$ $L^{q_n}$-weakly converges to $f$. We obtain the weak compactness in our setting. The proof is completed.
\end{proof} 

\begin{prop}\label{sumconv}
Let $f_n, g_n \in L^{q_n}(X_n, m_n)$ and $f, g \in L^q(Y, m)$. Assume that $f_n$, $g_n$ $L^{q_n}$-strongly converges to $f$, $g$ respectively. Then $f_n + g_n$ $L^{q_n}$-strongly converges to $f + g$ and it holds that
\begin{equation}\label{convnorm}
\lim_{n \to \infty} \| f_n \|_{L^{q_n}(X_n, m_n)} = \| f \|_{L^{q}(Y, m)}.
\end{equation}
\end{prop}

\begin{prop}\label{inproconv}
Let $f_n \in L^q(X_n, m_n)$, $f \in L^q(Y, m)$ and let $g_n \in L^{q^*}(X_n, m_n)$, $g \in L^{q^*}(Y, m)$ with $q^* = q/(q - 1)$. Assume that $f_n$ $L^q$-strongly converges to $f$ and $g_n$ $L^{q^*}$-weakly converges to $g$. Then it holds that
\begin{equation}\label{convprod}
\lim_{n \to \infty} \int_{X_n} f_n(x) g_n(x) \, dm_n(x) = \int_{Y} f(y) g(y) \, dm(y).
\end{equation}
\end{prop}

\begin{rem}
\begin{enumerate}
\item We do not need any continuity of $p_n$ for the definition of $L^{q_n}$-convergence and some properties as above. We only use the Borel measurability of $p_n$.

\item From the above properties, this $L^{q_n}$-convergence is an {\it asymptotic relation} which defined by \cite{KS}, so that this is regarded as a natural extension of $L^{q_n}$-convergence in \cite{GMS}*{Section 6.1} and \cite{AH}*{Section 3}. 
\end{enumerate}
\end{rem}

Under this $L^{q_n}$-convergence, we define convergences of $q$-Cheeger energy functionals.

\begin{dfn}\label{moscoCh}
We say that $\Ch_{q_n}^{X_n}$ $\Gamma$-converges to $\Ch_q^Y$ if 
\begin{enumerate}
\item for any sequence of functions $f_n \in L^{q_n}(X_n, m_n)$ $L^{q_n}$-strongly converging to a function $f \in L^q(Y, m)$, we have
\begin{equation*}
\liminf_{n \to \infty} \Ch_{q_n}^{X_n}(f_n) \geq \Ch_q^Y(f),
\end{equation*}
\item for any $f \in L^q (Y, m)$, there exists a sequence of functions $f_n \in L^{q_n}(X_n, m_n)$ $L^{q_n}$-strongly convergent to $f$ such that
\begin{equation*}
\lim_{n \to \infty} \Ch_{q_n}^{X_n}(f_n) = \Ch_q^Y(f).
\end{equation*}
\end{enumerate}
Moreover, we say that $\Ch_{q_n}^{X_n}$ Mosco converges to $\Ch_q^Y$ if
\begin{enumerate}
\item[{\rm (1$'$)}] for any sequence of functions $f_n \in L^{q_n}(X_n, m_n)$ $L^{q_n}$-weakly converging to a function $f \in L^q(Y, m)$, we have
\begin{equation*}
\liminf_{n \to \infty}  \Ch_{q_n}^{X_n}(f_n) \geq \Ch_q^Y(f),
\end{equation*}
\end{enumerate}
and above (2).
\end{dfn}

The following results obtained in \cite{GMS} and \cite{AH} mean the pmG-convergent case of Theorem \ref{intro_main} (2) and (4). We need it for the proof of Theorem \ref{intro_main}.

\begin{thm}[\cite{GMS}*{Theorem 6.8}]\label{GMSCh}
Let $\{X_n\}_{n \in \mathbb{N}}$ be a sequence of p.m.m.~spaces satisfying $\CD(K, \infty)$ for a common number $K \in \mathbb{R}$ and $Y$ a p.m.m.~space. Assume that $X_n$ pmG-converges to $Y$ as $n \to \infty$. Then $\Ch_2^{X_n}$ Mosco converges to $\Ch_2^Y$ as $n \to \infty$.
\end{thm}

\begin{thm}[\cite{AH}*{Theorem 8.1}]\label{AHCh}
Let $\{X_n\}_{n \in \mathbb{N}}$ be a sequence of p.m.m.~spaces satisfying $\RCD(K, \infty)$ for a common number $K \in \mathbb{R}$ and $Y$ a p.m.m.~space. Assume that $X_n$ pmG-converges to $Y$ as $n \to \infty$. Then $\Ch_{q_n}^{X_n}$ $\Gamma$-converges to $\Ch_q^Y$ as $n \to \infty$.
\end{thm}

\begin{proof}[{\bf Proof of Theorem \ref{intro_main} (2) and (4)}]

Let $\{q_n\} \subset (1, \infty)$ be a sequence convergent to a real number $q \in (1, \infty)$. We take any functions $f_n \in L^{q_n}(X_n, m_n), f \in L^q(Y, m)$ and assume that $f_n$ $L^{q_n}$-weakly converges to $f$ as $n \to \infty$. For each $n$, we define a function $g_n \in L^{q_n}(X_n^*, m_n^*)$ by
\begin{equation*}
g_n(y) := \int_{X_n} f_n(x) \, d\mu_y^n(x),
\end{equation*}
where $\{ \mu_y^n \}_{y \in X_n^*}$ is the disintegration of $m_n$ for $p_n$. By Lemma \ref{compare_wug},
\begin{equation*}
|Dg_n|_w(y)^{q_n} \leq \int_{X_n} |Df_n|_w(x)^{q_n} \, d\mu_y^n (x) 
\end{equation*}
holds for $m_n^*$-a.e.~$y \in X_n^*$. Thus we have
\begin{align}
\Ch_{q_n}^{X_n^*}(g_n) & = \, \frac{1}{q_n} \int_{X_n^*} |Dg_n|_w(y)^{q_n} \, dm_n^*(y) \nonumber \\
& \leq \, \frac{1}{q_n} \int_{X_n^*} \int_{X_n} |Df_n|_w(x)^{q_n} \, d\mu_y^n (x) dm_n^*(y) \label{ChYnChn} \\
& = \, \frac{1}{q_n} \int_{X_n} |Df_n|_w(x)^{q_n} \, dm_n(x) = \Ch_{q_n}^{X_n}(f_n). \nonumber
\end{align}
Since $f_n$ $L^{q_n}$-weakly converges to $f$, the sequence $\{g_n\}$ $L^q$-weakly converges to $f$ on $Z$ in the same way as in the proof of Proposition \ref{weaklscnorm}. 

In the case of $q_n = q = 2$, by Theorem \ref{GMSCh} under the $\CD(K, \infty)$ assumption, $\Ch_2^{X_n^*}$ Mosco converges to $\Ch_2^Y$. Combining this with (\ref{ChYnChn}) implies
\begin{equation*}
\Ch_2^Y(f) \leq \liminf_{n \to \infty} \Ch_2^{X_n^*}(g_n) \leq \liminf_{n \to \infty} \Ch_2^{X_n}(f_n).
\end{equation*}
We obtain ($1'$) in Definition \ref{moscoCh} in the case of $q_n = q = 2$.

In the general case, we further assume that $f_n$ $L^{q_n}$-strongly converges to $f$. We prove the following claim.

\begin{claim}\label{strconvclaim}
$g_n$ $L^{q_n}$-strongly converges to $f$ on $Z$ as $n \to \infty$.
\end{claim}

\begin{proof}
We take any $\varepsilon > 0$. Since $f_n$ $L^{q_n}$-strongly converges to $f$, there exists a decomposition $f_n = \alpha_n + \beta_n$ such that $\alpha_n$ $L^1$-strongly convergent and
\begin{equation*}
\sup_{n \in \N} \| \alpha_n \|_{L^\infty(X_n, m_n)} < +\infty, \quad \sup_{n \in \N} \| \beta_n \|_{L^{p_n}(X_n, m_n)} < \varepsilon.
\end{equation*}
We define $\hat{\alpha}_n, \hat{\beta}_n$ by 
\begin{equation*}
\hat{\alpha}_n(y) := \int_{X_n} \alpha_n(x) \, d\mu_y^n(x), \quad \hat{\beta}_n(y) := \int_{X_n} \beta_n(x) \, d\mu_y^n(x),
\end{equation*}
where $\{ \mu_y^n \}_{y \in X_n^*}$ is the disintegration of $m_n$ for $p_n$. We see that $g_n = \hat{\alpha}_n + \hat{\beta}_n$ and
\begin{equation*}
\| \hat{\alpha}_n \|_{L^\infty(X_n^*, m_n^*)} \leq \| \alpha_n \|_{L^\infty(X_n, m_n)}, \quad  \| \hat{\beta}_n \|_{L^{p_n}(X_n^*, m_n^*)} \leq \| \beta_n \|_{L^{p_n}(X_n, m_n)}.
\end{equation*}
Thus it suffices to prove the $L^1$-strongly convergence of $\{\hat{\alpha}_n\}$ on $Z$. Splitting $\alpha_n$ in the positive and negative parts, we assume $\alpha_n \geq 0$. Let $\alpha$ be the $L^1$-strong limit of $\{\alpha_n\}$ and let $\hat{\alpha}$ be a $L^2$-weak limit on $Z$ of a convergent subsequence of $\{\sigma \circ \hat{\alpha}_n\}$, where $\sigma(t) := \mathrm{sign}(t)\sqrt{|t|}$ is the signed root. For any $\varphi \in \Cbs(Z)$, we see that
\begin{align*}
\int_{X_n} \varphi(p_n(x)) (\sigma \circ \alpha_n)(x) \, dm_n(x) &= \int_{X_n^*} \varphi(y) \int_{X_n} \alpha_n(x)^{\frac{1}{2}} \, d\mu_y(x)dm_n^*(y) \\
&\leq \int_{X_n^*} \varphi(y) (\sigma \circ \hat{\alpha}_n)(y) \, dm_n^*(y),
\end{align*}
which implies, by taking the limit as $n \to \infty$,
\begin{equation*}
\int_{Y} \varphi(y) (\sigma \circ \alpha)(y) \, dm(y) \leq \int_{Y} \varphi(y) \hat{\alpha}(y) \, dm(y).
\end{equation*}
Thus we have $\sigma \circ \alpha(y) \leq \hat{\alpha}(y)$ for $m$-a.e.~$y \in Y$. On the other hand, we have
\begin{equation*}
\begin{split}
\|\hat{\alpha}\|_{L^2(Y, m)} &\leq \liminf_{n \to \infty} \| \sigma \circ \hat{\alpha}_n\|_{L^2(X_n^*, m_n^*)}\\
&\leq \liminf_{n \to \infty} \| \sigma \circ \alpha_n\|_{L^2(X_n, m_n)} = \|\sigma\circ\alpha\|_{L^2(Y, m)}.
\end{split}
\end{equation*}
Combining these leads to $\sigma \circ \alpha(y) = \hat{\alpha}(y)$ for $m$-a.e.~$y \in Y$ and
\begin{equation*}
\lim_{n \to \infty} \| \sigma \circ \hat{\alpha}_n\|_{L^2(X_n^*, m_n^*)} = \|\sigma\circ\alpha\|_{L^2(Y, m)}.
\end{equation*}
Therefore $\sigma \circ \hat{\alpha}_n$ $L^2$-strongly converges to $\sigma \circ \alpha$, that is, $\hat{\alpha}_n$ $L^1$-strongly converges to $\alpha$. This completes the proof.
\end{proof}

By Theorem \ref{AHCh} under the $\RCD(K, \infty)$ assumption, $\Ch_{q_n}^{X_n^*}$ $\Gamma$-converges to $\Ch_q^Y$. Combining this with Claim \ref{strconvclaim} and (\ref{ChYnChn}) implies
\begin{equation*}
\Ch_q^Y(f) \leq \liminf_{n \to \infty} \Ch_{q_n}^{X_n^*}(g_n) \leq \liminf_{n \to \infty} \Ch_{q_n}^{X_n}(f_n).
\end{equation*}
We obtain (1) in Definition \ref{moscoCh} in the general case.

We next prove (2) in Definition \ref{moscoCh}. We take any function $f \in L^q(Y, m)$. By Theorem \ref{GMSCh} and \ref{AHCh}, the sequence of the Cheeger energy functionals $\Ch_{q_n}^{X_n^*}$  $\Gamma$-converges to $\Ch_q^Y$ as $n \to \infty$. Thus, there exists a sequence of functions $g_n \in L^{q_n}(X_n^*, m_n^*)$ $L^{q_n}$-strongly convergent to $f$ on $Z$ such that
\begin{equation*}
\lim_{n \to \infty} \Ch_{q_n}^{X_n^*}(g_n) = \Ch_q^Y(f).
\end{equation*}
We define a function $f_n$ by $f_n := {p_n}^\ast g_n$ for each $n$. These $f_n$ are what we want. Let us prove that these $f_n$ satisfy the condition (2) in Definition \ref{moscoCh}. We see that
\begin{equation*}
\int_{X_n} f_n(x)^{q_n} \, dm_n(x) = \int_{X_n} g_n(p_n(x))^{q_n} \, dm_n(x) = \int_{X_n^*} g_n(y)^{q_n} \, dm_n^*(y)
\end{equation*}
holds, so that $f_n \in L^{q_n}(X_n, m_n)$. Moreover, we obtain
\begin{equation}\label{proof(2)_normconv}
\lim_{n \to \infty} \| f_n \|_{L^{q_n}(X_n, m_n)} = \| f \|_{L^{q}(Y, m)}.
\end{equation}
For any $\varphi \in \Cbs(Z)$, we have
\begin{align*}
\int_{X_n} \varphi(p_n(x))f_n(x) \, dm_n(x) = & \int_{X_n} \varphi(p_n(x)) g_n(p_n(x)) \, dm_n(x) \\
= & \int_{X_n^*} \varphi(y) g_n(y) \, dm_n^*(y).
\end{align*}
Thus we obtain
\begin{equation}\label{proof(2)_weakconv}
\lim_{n \to \infty} \int_{X_n} \varphi(p_n(x)) f_n(x) \, dm_n(x) = \int_Y \varphi(y) f(y) \, dm(y)
\end{equation}
for any $\varphi \in \Cbs(Z)$. Combining (\ref{proof(2)_normconv}) and (\ref{proof(2)_weakconv}) implies that $f_n$ $L^2$-strongly converges to $f$. By Theorem \ref{stabCh}, we have
\begin{equation*}
\Ch_{q_n}^{X_n}(f_n) = \Ch_{q_n}^{X_n^*}(g_n)
\end{equation*}
and so
\begin{equation*}
\lim_{n \to \infty} \Ch_{q_n}^{X_n}(f_n) = \Ch_q^Y(f).
\end{equation*}
Therefore the sequence of functions $f_n$ satisfies the conditions (2). This completes the proof of Theorem \ref{intro_main}.
\end{proof}

\subsection{Semicontinuity of spectra and spectral gaps}

As an application of the Mosco convergence of the Cheeger energy functionals, we obtain the semicontinuity of the spectra of Laplacians on p.m.m. spaces satisfying $\RCD(K, \infty)$. The Laplacian $\Delta_X$ on a p.m.m. space $X$ satisfying $\RCD(K, \infty)$ is defined as the self-adjoint linear operator associated with the quadratic form $\Ch_2^X$. We denote by $\sigma(\Delta_X)$ the spectrum of $\Delta_X$.

\begin{proof}[{\bf Proof of Corollary \ref{intro_spec}} ]
This corollary follows from Theorem \ref{intro_main} and \cite{KS}*{Proposition 5.30} directly.
\end{proof}

\begin{proof}[{\bf Proof of Corollary \ref{intro_spec_fc}}]
We take any $\lambda \in \sigma(\Delta_{X^*})$. By Corollary \ref{intro_spec}, there exists a sequence $\lambda_n \in \sigma(\Delta_X)$ convergent to $\lambda$. Since $\sigma(\Delta_X)$ is closed in $[0 , \infty)$, $\lambda$ belongs to $\sigma(\Delta_X)$. The proof is completed.
\end{proof}

Let $(X, d, m)$ be an m.m.~space with finite mass. For any real number $q \in (1, \infty)$ and any $f \in L^q(X, m)$, we set
\begin{equation}
c_q(f) := \left( \inf_{a \in \R} \int_X \left|f(x) - a\right|^q \, dm(x)\right)^{\frac{1}{q}}.
\end{equation}
We define the {\it $q$-spectral gap} $\lambda_{1, q}(X, d, m)$ by
\begin{equation}
\lambda_{1, q}(X, d, m) :=  \inf_{f} \frac{q\Ch_q^X(f)}{c_q(f)^q} ,
\end{equation}
where $f$ runs over all nonconstant $L^q$-functions. It is well-known that the infimum does not change if we minimize $q\Ch_q^X(f)$ where $f$ runs over all nonconstant $L^q$-functions with
\begin{equation*}
\| f \|_{L^q(X, m)} = 1, \quad \int_X \left|f(x)\right|^{q-2}f(x) \, dm(x) = 0.
\end{equation*}

\begin{prop}
Let $(X, d, m)$ be an m.m.~space with a metric measure foliation and let $q \in (1, \infty)$ be a real number. Assume that $m$ has the finite mass. Then we have
\begin{equation}
\lambda_{1, q}(X^*, d^*, m^*) \geq \lambda_{1, q}(X, d, m).
\end{equation}
\end{prop}

\begin{proof}
We take any $f \in L^q(X^*, m^*)$. By Theorem \ref{stabCh}, we have
\begin{equation*}
\Ch_q^{X^*}(f) = \Ch_q^X(p^* f).
\end{equation*}
Moreover, for any $a \in \R$, we see that
\begin{equation*}
\int_{X^*} \left|f(y) - a\right|^q \, dm^*(y) = \int_X \left|p^* f(x) - a\right|^q \, dm(x),
\end{equation*}
which implies that $c_q(f) = c_q(p^* f)$. Therefore we obtain
\begin{equation*}
\lambda_{1, q}(X, d, m) \leq \frac{q\Ch_q^X(p^* f)}{c_q(p^* f)^q} = \frac{q\Ch_q^{X^*}(f)}{c_q(f)^q}.
\end{equation*}
This completes the proof.
\end{proof}

\begin{prop}
Let $\{(X_n, d_n, m_n, \bar{x}_n)\}_{n \in \N}$ be a sequence of p.m.m. spaces and $(Y, d, m, \bar{y})$ be a p.m.m.~space. Let $\{q_n\}_{n \in \N} \subset (1, \infty)$ be a sequence of real numbers convergent to $q \in (1, \infty)$. Assume that the same assumptions as in Theorem \ref{intro_main} (3) and each $m_n$ has the finite mass. Then we have
\begin{equation}\label{limsupspec}
\limsup_{n \to \infty} \lambda_{1, q_n}(X_n, d_n, m_n) \leq \lambda_{1, q}(Y, d, m).
\end{equation}
\end{prop}

\begin{proof}
We take any $f \in L^q(Y, m)$. By Theorem \ref{intro_main} (4), there exists $f_n \in L^{q_n}(X_n, m_n)$ $L^{q_n}$-strongly converging to $f$ such that
\begin{equation*}
\lim_{n \to \infty} \Ch_{q_n}^{X_n}(f_n) = \Ch_q^Y(f).
\end{equation*}
Then we see that $c_{q_n}(f_n)$ converges to $c_{q}(f)$ in the same way as in \cite{AH}*{Lemma 9.2}. Thus we obtain
\begin{equation*}
\limsup_{n \to \infty} \lambda_{1, q_n}(X_n, d_n, m_n) \leq \limsup_{n \to \infty} \frac{q_n \Ch_{q_n}^{X_n}(f_n)}{c_{q_n}(f_n)^{q_n}} = \frac{q \Ch_{q}^{Y}(f)}{c_{q}(f)^q},
\end{equation*}
which leads to (\ref{limsupspec}).
\end{proof}

\begin{rem}
In our setting, we do not obtain the upper semicontinuity of the spectra $\sigma(\Delta_{X_n})$ and the lower semicontinuity of the $q_n$-spectral gap $\lambda_{1, q_n}(X_n, d_n, m_n)$. Via the metric measure foliation, the spectral information is lost in general. The limit-like space $Y$ does not have the full informations of the limit behavior of $X_n$.
\end{rem}

\begin{ex}
Consider the sequence $\{S^n \times S^1\}_{n \in \N}$ of the Riemannian product of the $n$-dimensional unit sphere $S^n$ and the $1$-dimensional unit sphere $S^1$. The product space $S^n \times S^1$ has a metric measure foliation induced by the trival $S^1$-fibration. We see that
\begin{equation*}
\lambda_{1, 2}(S^n) = n, \quad \lambda_{1,2}(S^n \times S^1) = \lambda_{1,2}(S^1) = 1.
\end{equation*}
\end{ex}

\begin{ex}
Consider the sequence $\{S^n(\sqrt{n-1})\}_{n \in \N}$, the map $p_n$ of (\ref{sphere_mmf}), and the space $X_n$ in Section 1. $X_n$ pmG-converges to 1-dimensional standard Gaussian space $(\R, |\cdot|, \gamma)$ as $n \to \infty$ (see \cite{Nkj}*{Lemma 3.9}), where $\gamma$ is the 1-dimensional standard Gaussian measure (i.e.~$\gamma = \gamma_1$). It is well-known that
\begin{equation*}
\begin{split}
\sigma(\Delta_{S^{n}(\sqrt{n-1})}) &= \left\{ k \left(1 + \frac{k}{n-1} \right) \midd  k = 0, 1, 2, \ldots \right\}, \\
\sigma(\Delta_{(\R, |\cdot|, \gamma)})& = \left\{ k \, | \, k = 0, 1, 2, \ldots \right\}.
\end{split}
\end{equation*}
By Corollary \ref{intro_spec_fc}, we see that $\sigma(\Delta_{X_n}) \subset \sigma(\Delta_{S^{n}(\sqrt{n-1})})$. Moreover, by \cite{GMS}*{Theorem 7.8}, $k$-th eigenvalue of $\Delta_{X_n}$ converges to $k$-th eigenvalue of $\Delta_{(\R, |\cdot|, \gamma)}$ taking account of the multiplicity. Thus we see that $\sigma(\Delta_{X_n}) = \sigma(\Delta_{S^{n}(\sqrt{n-1})})$ and each multiplicity of eigenvalues of $\Delta_{X_n}$ equals to 1.
\end{ex}

\subsection{Descending slope of the relative entropy and heat flow}

In this subsection, we state the results for the (descending) slope of the relative entropy and the heat flow. We first define these notions along \cites{calc, GMS}. Let $(X, d, m)$ be an m.m.~space.

\begin{dfn}[Descending slope of $\Ent_m$]\label{DEnt}
We define the ({\it descending}) {\it slope} $\DEnt{m} : \Pb_2(X) \to [0, +\infty]$ of the relative entropy $\Ent_m$ for $m$ by
\begin{equation}
\DEnt{m}(\mu) := \limsup_{W_2(\nu, \mu) \to 0} \frac{(\Ent_m(\mu) - \Ent_m(\nu))^+}{W_2(\mu, \nu)},
\end{equation}
where $(\cdot)^+$ means the positive part. $\DEnt{m}(\mu)$ is equal to $+\infty$ if $\mu$ is outside $D(\Ent_m)$ and $0$ if $\mu$ is an isolated measure in $D(\Ent_m)$.
\end{dfn}

\begin{prop}\label{CDDEnt}
Let $(X, d, m)$ be an m.m.~space satisfying $\CD(K, \infty)$ for a real number $K \in \mathbb{R}$. Then, for any $\mu \in \Pb_2(X)$, we have
\begin{equation}
\DEnt{m}(\mu) = \sup_{\mu \neq \nu \in \Pb_2(X)} \left(\frac{\Ent_m(\mu) - \Ent_m(\nu)}{W_2(\mu, \nu)} + \frac{K}{2} W_2(\mu, \nu) \right)^+.
\end{equation}
In particular, $\DEnt{m}$ is lower semicontinuous with respect to $W_2$.
\end{prop}

It is known that the slope $\DEnt{m}$ of $\Ent_m$ and the 2-Cheeger energy functional $\Ch_2$ on an m.m.~space satisfying $\CD(K, \infty)$ have the following relation. This relation is a deep result obtained in \cite{calc}.

\begin{thm}[\cite{calc}*{Theorem 7.6}]\label{chslope}
Let $(X, d, m)$ be an m.m.~space satisfying $\CD(K, \infty)$ for a number $K \in \mathbb{R}$. Then, for any measure $\mu \in \Pb_2(X)$ that is absolutely continuous with respect to $m$, we have 
\begin{equation}
\DEnt{m}^2(\mu) = 8\Ch_2(\sqrt{\rho}),
\end{equation}
where $\rho$ is the density of $\mu$ with respect to $m$.
\end{thm}

Let $I \subset \mathbb{R}$ be an interval of $\mathbb{R}$. A curve $\gamma : I \to X$ is said to be {\it locally absolutely continuous} if there exists $f \in L^1_{loc}(I)$ satisfying (\ref{abs_conti}) for any $s, t \in I$ with $s < t$. We denote by $AC_{loc}(I; X)$ the set of locally absolutely continuous curves on $X$. For each $\gamma \in AC_{loc}(I; X)$, the metric derivative $|\dot{\gamma}| \in L^1_{loc}(I)$ of $\gamma$ is defined by (\ref{metderi}) locally and is the minimal function, in the a.e. sense, of $L^1_{loc}$-functions satisfying (\ref{abs_conti}). 

\begin{dfn}[Gradient flow of the relative entropy]\label{gradflow}
Let $(X, d, m)$ be an m.m.~space satisfying $\CD(K, \infty)$ for a number $K \in \mathbb{R}$ and let $\bar{\mu} \in \Pb_2(X) \cap D(\Ent_m)$. A curve $\mu : [0, +\infty) \to \Pb_2(X) \cap D(\Ent_m)$ is a {\it $W_2$-gradient flow of $\Ent_m$ starting from $\bar{\mu}$} provided $\mu$ belongs to $AC_{loc}\left([0, +\infty); (\Pb_2(X), W_2)\right)$ and satisfies $\mu_0 = \overline{\mu}$ and 
\begin{equation}\label{EDE}
\Ent_m(\mu_s) = \Ent_m(\mu_t) + \frac{1}{2} \int_s^t |\dot{\mu}|(r)^2 \, dr + \frac{1}{2} \int_s^t \DEnt{m}(\mu_r)^2 \, dr
\end{equation} 
for any $s, t \in [0, +\infty)$ with $s < t$.
\end{dfn}

\begin{rem}
The formula (\ref{EDE}) is called the {\it Energy Dissipation Equation} (abbreviated as EDE) and a gradient flow in Definition \ref{gradflow} is called a {\it flow in the EDE sense}. Moreover a $W_2$-gradient flow $\mu$ of $\Ent_m$ satisfies 
\begin{equation*}
-\frac{d}{dt} \Ent_m(\mu_t) = \DEnt{m}(\mu_t)^2
\end{equation*}
for a.e.~$t \in (0, \infty)$. As one of the most important results in \cite{calc}, it is known that the $W_2$-gradient flows of $\Ent_m$ coincide with the $L^2$-gradient flows of the 2-Cheeger energy functional $\Ch_2$ on an m.m.~space satisfying $\CD(K, \infty)$. Thus $\mu$ is called a {\it heat flow}.
\end{rem}

The existence and uniqueness of the $W_2$-gradient flow of $\Ent_m$ was proved by \cite{Ggradflow} in the case of locally compact space and by \cite{calc} in the general case.

\begin{thm}[\cite{calc}*{Theorem 9.3}]
Let $(X, d, m)$ be an m.m.~space satisfying $\CD(K, \infty)$ for a number $K \in \mathbb{R}$. Then, for any measure $\bar{\mu} \in \Pb_2(X) \cap D(\Ent_m)$, there exists a unique $W_2$-gradient flow of $\Ent_m$ starting from $\bar{\mu}$.
\end{thm}

In order to describe Theorem \ref{intro_main2} more precisely, we define a $\Gamma$-convergence of the slopes of the relative entropies in our setting.

\begin{dfn} \label{convDEnt}
Under the assumption of Theorem \ref{intro_main2}, we say that $\DEnt{m_n}$ $\Gamma$-converges to $\DEnt{m}$ if
\begin{enumerate}
\item For any $\mu_n \in \Pb_2(X_n), \mu \in \Pb_2(Y)$ such that $W_2(\mu_n, {p_n}^* \mu)$ tends to $0$ as $n \to \infty$, we have 
\begin{equation}
\liminf_{n \to \infty} \DEnt{m_n}(\mu_n) \geq \DEnt{m}(\mu).
\end{equation}
\item For any $\mu \in \Pb_2 (Y)$, there exists a sequence $\mu_n \in \Pb_2(X_n)$ such that
\begin{align}
& \lim_{n \to \infty} W_2(\mu_n, {p_n}^* \mu) = 0, \\
& \lim_{n \to \infty} \DEnt{m_n}(\mu_n) = \DEnt{m}(\mu).
\end{align}
\end{enumerate}
\end{dfn}

\begin{rem}
In Definition \ref{convDEnt},  the definition of the convergence of $\mu_n \in \Pb_2(X_n)$ to $\mu \in \Pb_2(Y)$ is regarded as
\begin{equation}\label{conv_meas_pull}
\lim_{n \to \infty} W_2(\mu_n, {p_n}^* \mu) = 0.
\end{equation}
This convergence is an asymptotic relation in \cite{KS}, so that this is a natural convergence. On the other hand, the $W_2$-convergence of the push-forward measures ${p_n}_* \mu_n$ is weaker than (\ref{conv_meas_pull}) and is not an asymptotic relation.
\end{rem}

We need the following lemma for the proof of Theorem \ref{intro_main2}.

\begin{lem}\label{stabDEnt}
Let $(X, d, m)$ be an m.m.~space with a metric measure foliation. Then, for any $\mu \in \Pb_2(X^*)$, we have
\begin{equation}\label{stabDEnt_eq}
\DEnt{m}(p^* \mu) = \DEnt{m^*}(\mu).
\end{equation}
\end{lem}

\begin{proof}
We take any measure $\mu \in \Pb_2(X^*)$. If $\mu$ is outside $D(\Ent_{m^*})$ or is an isolated measure in $D(\Ent_{m^*})$, the pullback measure $p^* \mu$ is also the same, which implies (\ref{stabDEnt_eq}) trivially. We only prove the case that $\mu$ belongs to $D(\Ent_{m^*})$ and is non-isolated. By (\ref{pushent}), Proposition \ref{prop_pull} (3), and Lemma \ref{pushW_q}, we have
\begin{align*}
\DEnt{m}(p^* \mu) & = \limsup_{W_2(\nu, p^* \mu) \to 0} \frac{(\Ent_{m}(p^* \mu) - \Ent_{m}(\nu))^+}{W_2(p^* \mu, \nu)} \\
& \leq \limsup_{W_2(\nu, p^* \mu) \to 0} \frac{(\Ent_{m^*}(\mu) - \Ent_{m^*}(p_* \nu))^+}{W_2(\mu, p_* \nu)} \\
& \leq \limsup_{W_2(\nu', \mu) \to 0} \frac{(\Ent_{m^*}(\mu) - \Ent_{m^*}(\nu'))^+}{W_2(\mu, \nu')} \\
& = \DEnt{m^*}(\mu).
\end{align*}
On the other hand, by Proposition \ref{prop_pull} (3) and Lemma \ref{pullW_q}, we have
\begin{align*}
\DEnt{m}(p^* \mu) & = \limsup_{W_2(\nu, p^* \mu) \to 0} \frac{(\Ent_{m}(p^* \mu) - \Ent_{m}(\nu))^+}{W_2(p^* \mu, \nu)} \\
& \geq \limsup_{W_2(\nu', \mu) \to 0} \frac{(\Ent_{m}(p^* \mu) - \Ent_{m}(p^* \nu'))^+}{W_2(p^* \mu, p^* \nu')} \\
& = \limsup_{W_2(\nu', \mu) \to 0} \frac{(\Ent_{m^*}(\mu) - \Ent_{m^*}(\nu'))^+}{W_2(\mu, \nu')} \\
& = \DEnt{m^*}(\mu).
\end{align*}
Therefore we obtain (\ref{stabDEnt_eq}). This completes the proof.
\end{proof}

Moreover, we need the following two results obtained in \cite{GMS}.

\begin{thm}[\cite{GMS}*{Theorem 4.7}]\label{GMSEnt}
Let $\{(X_n, d_n, m_n, \bar{x}_n)\}_{n \in \mathbb{N}}$ be a sequence of p.m.m.~spaces and $(Y, d, m, \bar{y})$ a p.m.m.~space. Assume that $X_n$ pmG-converges to $Y$ as $n \to \infty$. Then, $\Ent_{m_n}$ $\Gamma$-converges to $\Ent_m$ as $n \to \infty$, that is, the following (1) and (2) hold.
\begin{enumerate}
\item For any sequence of measures $\mu_n \in \Pb_2(X_n)$ $W_2$-converging to a measure $\mu \in \Pb_2(Y)$, we have 
\begin{equation*}
\liminf_{n \to \infty} \Ent_{m_n}(\mu_n) \geq \Ent_{m}(\mu).
\end{equation*}
\item For any $\mu \in \Pb_2 (Y)$, there exists a sequnece of measures $\mu_n \in \Pb_2(X_n)$ $W_2$-convergent to $\mu$ such that
\begin{equation*}
\lim_{n \to \infty} \Ent_{m_n}(\mu_n) = \Ent_{m}(\mu).
\end{equation*}
\end{enumerate}
Note that the $W_2$-convergence of $\mu_n$ is well-defined since each $\mu_n$ and $\mu$ are regarded as measures on the common $Z$.
\end{thm}

\begin{thm}[\cite{GMS}*{Theorem 5.14}]\label{GMSDEnt}
Let $\{(X_n, d_n, m_n, \bar{x}_n)\}_{n \in \mathbb{N}}$ be a seqeuence of p.m.m.~spaces satisfying $\CD(K, \infty)$ for a common number $K \in \mathbb{R}$ and $(Y, d, m, \bar{y})$ a p.m.m.~space. Assume that $X_n$ pmG-converges to $Y$ as $n \to \infty$. Then $\DEnt{m_n}$ Mosco converges to $\DEnt{m}$ as $n \to \infty$, that is, the following (1) and (2) hold.
\begin{enumerate}
\item For any sequence of measures $\mu_n \in \Pb_2(X_n)$ weakly converging to a measure $\mu \in \Pb_2(Y)$, we have 
\begin{equation*}
\liminf_{n \to \infty} \DEnt{m_n}(\mu_n) \geq \DEnt{m}(\mu).
\end{equation*}
\item For any $\mu \in \Pb_2 (Y)$, there exists a sequnece of measures $\mu_n \in \Pb_2(X_n)$ $W_2$-convergent to $\mu$ such that
\begin{equation*}
\lim_{n \to \infty} \DEnt{m_n}(\mu_n) = \DEnt{m}(\mu).
\end{equation*}
\end{enumerate}
Note that the weak and $W_2$ convergences of $\mu_n$ are well-defined since each $\mu_n$ and $\mu$ are regarded as measures on the common $Z$.
\end{thm}

\begin{proof}[{\bf Proof of Theorem \ref{intro_main2}}]
We first prove (1) in Definition \ref{convDEnt}. We take any measures $\mu_n \in \Pb_2(X_n), \mu \in \Pb_2(Y)$ and assume that
\begin{equation*}
\lim_{n \to \infty} W_2(\mu_n, {p_n}^* \mu) = 0.
\end{equation*}
We take any measure $\nu \in D(\Ent_{m})$. By Theorem \ref{GMSEnt} (2), there exists a sequence $\{ \nu_n \}_{n \in \mathbb{N}} \subset \Pb_2(X_n^*)$ $W_2$-convergent to $\nu$ such that $\Ent_{m_n^*}(\nu_n) \to \Ent_{m}(\nu)$. Then we have $W_2(\mu_n, {p_n}^* \nu_n) \to W_2(\mu, \nu)$ as $n \to \infty$ and, by Theorem \ref{GMSEnt} (1),
\begin{equation*}
\Ent_m(\mu) \leq \liminf_{n \to \infty} \Ent_{m_n^*}({p_n}_* \mu_n) \leq \liminf_{n \to \infty} \Ent_{m_n}(\mu_n).
\end{equation*}
Thus, by Proposition \ref{CDDEnt}, we have
\begin{align*}
& \frac{\Ent_m(\mu) - \Ent_m(\nu)}{W_2(\mu, \nu)} + \frac{K}{2}W_2(\mu, \nu) \\
& \leq \liminf_{n \to \infty} \left(\frac{\Ent_{m_n}(\mu_n) - \Ent_m({p_n}^* \nu_n)}{W_2(\mu_n, {p_n}^* \nu_n)} + \frac{K}{2}W_2(\mu_n, {p_n}^* \nu_n) \right) \\
& \leq \liminf_{n \to \infty} \DEnt{m_n}(\mu_n),
\end{align*}
which implies
\begin{equation*}
\liminf_{n \to \infty} \DEnt{m_n}(\mu_n) \geq \DEnt{m}(\mu).
\end{equation*}
The proof of (1) in Definition \ref{convDEnt} is completed.

We next prove (2) in Definition \ref{convDEnt}. We take any $\mu \in \Pb_2 (Y)$. By Theorem \ref{GMSDEnt} (2), there exists a sequence $\{ \nu_n \}_{n \in \mathbb{N}} \subset \Pb_2(X_n^*)$ $W_2$-convergent to $\nu$ such that $\DEnt{m_n^*}(\nu_n) \to \DEnt{m}(\mu)$ as $n \to \infty$. We define a measure $\mu_n$ by $\mu_n := {p_n}^* \nu_n \in \Pb_2(X_n)$ for each $n \in \mathbb{N}$. Then we have
\begin{equation*}
W_2(\mu_n, {p_n}^* \mu) = W_2(\nu_n, \mu) \to 0
\end{equation*}
and, by Lemma \ref{stabDEnt}, 
\begin{equation*}
\DEnt{m_n}(\mu_n) = \DEnt{m_n^*}(\nu_n) \to \DEnt{m}(\mu)
\end{equation*}
as $n \to \infty$. This completes the proof of (2) in Definition \ref{convDEnt}. We obtain Theorem \ref{intro_main2}.
\end{proof}

\begin{rem}
As in Theorem \ref{GMSDEnt}, the Mosco convergence of the slopes $\DEnt{m_n}$ is obtained in the pmG-convergent case. However, we do not know if we can extend Theorem \ref{intro_main2} to a suitable Mosco convergence in our setting.
\end{rem}

We obtain the following results about the heat flow.

\begin{prop}
Let $(X, d, m)$ be an m.m.~space with a metric measure foliation. Assume that $X$ satisfies $\CD(K, \infty)$ for a real number $K \in \mathbb{R}$ and its quotient space $X^*$ satisfies $\VG$. Let $\bar{\mu} \in \Pb_2(X^*) \cap D(\Ent_{m^*})$ and let $\mu : [0, +\infty) \to \Pb_2(X^*) \cap D(\Ent_{m^*})$ be the heat flow starting from $\bar{\mu}$. We define a curve $p^* \mu : [0, +\infty) \to \Pb_2(X) \cap D(\Ent_{m})$ by $(p^* \mu)_t := p^* \mu_t$. Then $p^* \mu$ is the heat flow starting from $p^* \bar{\mu}$.
\end{prop}

\begin{proof}
By the definition of the heat flow, it suffices to prove that
\begin{align}
\Ent_{m}(p^* \mu_t)& = \Ent_{m^*}(\mu_t) & \text{ for all } t > 0, \label{heatflow1} \\
|\dot{(p^* \mu)}|(t) & = |\dot{\mu}|(t) & \text{ for a.e. } t > 0,  \label{heatflow2} \\
\DEnt{m}(p^* \mu_t) & = \DEnt{m^*}(\mu_t) & \text{ for a.e. } t > 0.  \label{heatflow3}
\end{align}
(\ref{heatflow1}) and (\ref{heatflow3}) have already been obtained by Proposition \ref{prop_pull} and Lemma \ref{stabDEnt}. we prove (\ref{heatflow2}). By Lemma \ref{pullW_q}, we have
\begin{equation*}
W_2(p^* \mu_s, p^* \mu_t) = W_2(\mu_s, \mu_t)
\end{equation*}
for any $s, t > 0$, which implies (\ref{heatflow2}). This completes the proof.
\end{proof}

The following lemma is a generalization of \cite{GMS}*{Theorem 5.7}. However, the proof of this lemma is exactly the same as that of \cite{GMS}*{Theorem 5.7}.

\begin{lem}\label{convhf}
Let $\{(X_n, d_n, m_n, \bar{x}_n)\}_{n \in \mathbb{N}}$ be a sequence of p.m.m.~spaces and $(Y, d, m, \bar{y})$ a p.m.m.~space and let $K \in \R$. Assume that each $X_n$ satisfies $\CD(K, \infty)$ and has a metric measure foliation and its quotient space $X_n^*$ satisfies $\VG$ and pmG-converges to $Y$ as $n \to \infty$. Let $\bar{\mu}_n \in \Pb_2(X_n) \cap D(\Ent_{m_n})$, $\bar{\mu} \in \Pb_2(Y) \cap D(\Ent_{m})$ and let $\mu_n, \mu$ be the heat flows starting from $\bar{\mu}_n, \bar{\mu}$ respectively and assume that  
\begin{equation*}
{p_n}_* \bar{\mu}_n \xrightarrow{W_2} \bar{\mu}, \quad \Ent_{m_n}(\bar{\mu}_n) \to \Ent_{m}(\bar{\mu}).
\end{equation*}
Then, for any $t > 0$, we have
\begin{equation}\label{push_conv_ heatflow}
{p_n}_* \mu_{n, t} \xrightarrow{W_2} \mu_t.
\end{equation}
\end{lem}

\begin{rem}
We conjecture that we could change the conclusion of Lemma \ref{convhf} from (\ref{push_conv_ heatflow}) to
\begin{equation*}
\lim_{n \to \infty} W_2(\mu_{n,t}, {p_n}^* \mu_t) = 0
\end{equation*}
if we assume the stronger convergence
\begin{equation*}
\lim_{n \to \infty} W_2(\bar{\mu}_n, {p_n}^* \bar{\mu}) = 0.
\end{equation*}
In the $\RCD(K, \infty)$ case, it is known that we have the contraction property
\begin{equation*}
W_2(\mu_{n,t}, {p_n}^* \mu_t) \leq e^{-Kt} W_2(\bar{\mu}_n, {p_n}^* \bar{\mu})
\end{equation*}
for any $t \geq 0$, so that this conjecture is true. However, we do not know if the conjecture is true in the general case.
\end{rem}

%% file: appendix.tex
\appendix
\section{Proof of Lemma \ref{CDlem}}

\begin{proof}[{\bf Proof of Lemma \ref{CDlem}}]
We take any $\mu_0, \mu_1 \in \Pb_V(X) \cap D(\Ent_m)$ with $W_2(\mu_0, \mu_1) < +\infty$ and any sufficiently small real number $\varepsilon > 0$. Let $\pi \in \Pb(X \times X)$ be an optimal transport plan for $W_2(\mu_0, \mu_1)$. By the tightness of $\pi$, there exist compact sets $K_0 \subset K_1 \subset \cdots \subset X \times X$ such that $\pi(K_0) \geq e^{-\varepsilon}$ and $\pi(K_n) \geq 1- \varepsilon e^{-n}$ for $n \geq 1$. Setting $A_0 := K_0$ and $A_n := K_n \setminus K_{n-1}$ for $n \geq 1$, we see that $\pi(A_n) \leq \varepsilon e^{-(n-1)}$ for $n \geq 1$ and
\begin{equation}\label{CDlemeq}
\theta(\varepsilon) := -\sum_{n = 0}^\infty \pi(A_n) \log{\pi(A_n)} \leq  \varepsilon + \sum_{n = 0}^\infty \varepsilon e^{-n} (\log{\varepsilon} + n) \to 0
\end{equation}
as $\varepsilon \to 0$. The inequality of (\ref{CDlemeq}) follows from the monotonicity of the function $r \mapsto r \log{r}$ for any sufficiently small $r$. We define the probability measures
\begin{equation*}
\pi_n := \pi(A_n)^{-1} \pi|_{A_n}, \quad \mu_0^n := {\pr_0}_* \pi_n, \quad \mu_1^n := {\pr_1}_* \pi_n
\end{equation*}
for $n \in \mathbb{N} \cup \{0\}$. Then we have $\mu_0^n, \mu_1^n \in \Pb_2(X) \cap D(\Ent_m)$ for each $n \in \mathbb{N} \cup \{0\}$ and
\begin{equation*}
\pi = \sum_{n=0}^\infty \pi(A_n)\pi_n, \quad \mu_0 = \sum_{n=0}^\infty \pi(A_n)\mu_0^n, \quad \mu_1 = \sum_{n=0}^\infty \pi(A_n)\mu_1^n.
\end{equation*}
We verify only $\mu_0^n, \mu_1^n \in D(\Ent_m)$. Setting $\rho_i, \rho_i^n$ the densities of $\mu_i, \mu_i^n$ respectively for $i = 0, 1$ and $n \in  \mathbb{N} \cup \{0\}$, we have $\rho_i^n(x) \leq \pi(A_n)^{-1}\rho_i(x)$ for $m$-a.e.~$x \in X$. Thus we have
\begin{align*}
& \int_{\{\rho_i^n > 1\}} \rho_i^n(x) \log{\rho_i^n(x)} \, dm(x) \\
& \leq \int_{\{\rho_i^n > 1\}} \rho_i^n(x) (\log{\rho_i(x)} - \log{\pi(A_n)}) \, dm(x) \\
& \leq \pi(A_n)^{-1} \int_{\{\rho_i^n(x_i) > 1\} \cap A_n} \log{\rho_i(x_i)} \, d\pi(x_0, x_1) - \log{\pi(A_n)} \\
& \leq \pi(A_n)^{-1} \int_{\{\rho_i(x_i) > 1\}} \log{\rho_i(x_i)} \, d\pi(x_0, x_1) - \log{\pi(A_n)} \\
& = \pi(A_n)^{-1} \int_{\{\rho_i > 1\}} \rho_i(x) \log{\rho_i(x)} \, d\mu_i(x) - \log{\pi(A_n)} < +\infty,
\end{align*}
which implies $\mu_0^n, \mu_1^n \in D(\Ent_m)$.
Therefore, by $\CD(K, \infty)$, there exists a $W_2$-geodesic $\mu^n : [0,1] \ni t \mapsto \mu_t^n \in \Pb_2(X)$ joining $\mu_0^n$ and $\mu_1^n$ satisfying (\ref{Kconvex}) for each $n \in \mathbb{N}\cup \{0\}$. We define a propability measure $\mu_t = \mu_t^\varepsilon$ for $t \in (0, 1)$ by
\begin{equation*}
\mu_t := \sum^\infty_{n=0} \pi(A_n) \mu_t^n.
\end{equation*}
It suffices to prove that $\mu : [0,1] \ni t \mapsto \mu_t \in \Pb(X)$ is a $W_2$-geodesic joining $\mu_0$ and $\mu_1$ satisfying 
\begin{equation}\label{CDlem_Kconvex}
\Ent_m (\mu_t) \leq (1-t) \Ent_m (\mu_0) + t \Ent_m (\mu_1) - \frac{K}{2} t(1-t) W_2(\mu_0, \mu_1)^2 + \theta(\varepsilon).
\end{equation}
In fact, by (\ref{CDlem_Kconvex}), we see that $\sup_{\varepsilon > 0} \Ent_m(\mu_t^\varepsilon) < +\infty$. Combining this and (\ref{PVW2}) and (\ref{Entformula}), we have $\sup_{\varepsilon > 0} \Ent_{\tilde{m}}(\mu_t^\varepsilon) < +\infty$. Thus $\{\mu_t^\varepsilon\}_{\varepsilon>0}$ is tight and then there exists a weak limit $\mu_t \in \Pb(X)$ of subsequence of $\{\mu_t^\varepsilon\}_{\varepsilon>0}$ as $\varepsilon \to 0$ for $t \in (0, 1)$. These weak limits $\{ \mu_t\}_{t \in (0, 1)}$ is also a $W_2$-geodesic and satisfies (\ref{Kconvex}) since $\theta(\varepsilon) \to 0$ as $\varepsilon \to 0$.

We first prove that $\mu$ is a $W_2$-geodesic. Since $\pi$ is an optimal transport plan, $\pi_n$ is also optimal for $W_2(\mu_0^n, \mu_1^n)$ for each $n \in \mathbb{N}\cup\{0\}$. Thus we have
\begin{align*}
\sum^\infty_{n=0} \pi(A_n) W_2(\mu_0^n, \mu_1^n)^2 &= \sum^\infty_{n=0} \pi(A_n) \int_{X\times X}d(x, x')^2 \, d\pi_n(x, x') \\
& = \sum^\infty_{n=0} \int_{A_n}d(x, x')^2 \, d\pi(x, x') = W_2(\mu_0, \mu_1)^2.
\end{align*}
By the triangle inequality, it is sufficient to prove that 
\begin{equation}\label{CDlem_geod}
W_2(\mu_t, \mu_i) \leq t^i (1-t)^{1-i} W_2(\mu_0, \mu_1)
\end{equation}
for $i=0, 1$. Let $\pi^{t,i}_n$ be an optimal transport plan for $W_2(\mu_t^n, \mu_i^n)$. Defining a measure
\begin{equation*}
\pi^{t,i} := \sum^\infty_{n=0} \pi(A_n) \pi^{t,i}_n,
\end{equation*}
we see that $\pi^{t,i} \in \Pi(\mu_t, \mu_i)$. Therefore, since $\mu^n$ is a $W_2$-geodesic,
\begin{align*}
W_2(\mu_t, \mu_i)^2 & \leq \int_{X\times X} d(x, x')^2 \, d\pi^{t,i} = \sum^\infty_{n=0} \pi(A_n) W_2(\mu_t^n, \mu_i^n)^2 \\
&= \sum^\infty_{n=0} \pi(A_n) t^{2i} (1-t)^{2(1-i)}W_2(\mu_0^n, \mu_1^n)^2 \\
& = t^{2i} (1-t)^{2(1-i)} W_2(\mu_0, \mu_1)^2,
\end{align*}
which implies (\ref{CDlem_geod}).

We next prove that $\mu$ satisfies (\ref{CDlem_Kconvex}). Let us first prove that
\begin{equation}\label{CDlem_ent}
\sum^\infty_{n=0} \pi(A_n) \Ent_m(\mu_i^n) \leq \Ent_m(\mu_i) + \theta(\varepsilon) < +\infty
\end{equation}
for $i= 0, 1$, where the series in the left-hand side converges. Since $\Ent_{\tilde{m}}(\mu_i^n) \geq 0$ and $\rho_i^n(x) \leq \pi(A_n)^{-1}\rho_i(x)$ for $m$-a.e.~$x \in X$, we have
\begin{align*}
(0 \leq) &\sum^\infty_{n=0} \pi(A_n) \Ent_{\tilde{m}}(\mu_i^n) = \sum^\infty_{n=0} \pi(A_n) \int_X \left( ze^{V^2}\rho_i^n\right) \log{\left(ze^{V^2}\rho_i^n \right)} \, d\tilde{m} \\
& \leq \sum^\infty_{n=0} \pi(A_n) \int_X \log{\left(ze^{V^2}\rho_i \right)} \, d\mu_i^n - \sum_{n = 0}^\infty \pi(A_n) \log{\pi(A_n)}\\
&= \sum^\infty_{n=0}  \int_{A_n} \log{\left(ze^{V^2}\rho_i \right)} \, d\pi + \theta(\varepsilon) = \int_X \log{\left(ze^{V^2}\rho_i \right)} \, d\mu_i + \theta(\varepsilon)\\
&=\Ent_m(\mu_i) + \int_X V^2 \, d\mu_i +\log{z}  + \theta(\varepsilon) < + \infty.
\end{align*}
Moreover, we obtain
\begin{equation*}
(0 \leq) \sum^\infty_{n=0} \pi(A_n) \int_X V^2 \, d\mu_i^n = \sum^\infty_{n=0} \int_{A_n} V^2 \, d\pi = \int_X V^2 \, d\mu_i < + \infty.
\end{equation*}
Therefore, the series
\begin{equation*}
\sum^\infty_{n=0} \pi(A_n) \Ent_m(\mu_i^n) = \sum^\infty_{n=0} \pi(A_n) \left( \Ent_{\tilde{m}}(\mu_i^n) - \int_X V^2 \, d\mu_i^n -\log{z} \right)
\end{equation*}
converges and (\ref{CDlem_ent}) holds. Let $\rho_t, \rho_t^n$ be the densities of $\mu_t, \mu_t^n$ respectively for $t \in (0, 1)$. By the definition of $\mu_t$ and Fubini's theorem, we have
\begin{equation*}
\rho_t(x) = \sum^\infty_{n=0} \pi(A_n) \rho_t^n(x)
\end{equation*}
for $m$-a.e.~$x \in X$. Therefore, by Jensen's inequality, Fubini's theorem, (\ref{Kconvex}) of $\mu^n$, and (\ref{CDlem_ent}),
\begin{align*}
&\Ent_m(\mu_t) \\
= & \int_X \left( \sum^\infty_{n=0} \pi(A_n) \rho_t^n(x) \right) \log{\left( \sum^\infty_{n=0} \pi(A_n) \rho_t^n(x) \right)} \, dm(x) \\
\leq & \int_X  \sum^\infty_{n=0} \pi(A_n) (\rho_t^n(x) \log{ \rho_t^n(x)} ) \, dm(x) = \sum^\infty_{n=0} \pi(A_n) \Ent_m(\mu_t^n) \\
\leq & (1-t) \sum^\infty_{n=0} \pi(A_n) \Ent_m (\mu_0^n) + t \sum^\infty_{n=0} \pi(A_n) \Ent_m (\mu_1^n) \\
& - \frac{K}{2} t(1 - t) \sum^\infty_{n=0} \pi(A_n) W_2(\mu_0^n, \mu_1^n)^2 \\
\leq & (1-t) \Ent_m (\mu_0) + t \Ent_m (\mu_1) - \frac{K}{2} t(1 - t) W_2(\mu_0, \mu_1)^2 + \theta(\varepsilon).
\end{align*}
We obtain (\ref{CDlem_Kconvex}). The proof is completed.
\end{proof}